\documentclass[11pt]{article}

\usepackage{amssymb}
\usepackage{amscd}
\usepackage{amsmath}
\usepackage{amsfonts}
\usepackage{amsthm}
\usepackage{color}
\usepackage{enumerate}
\usepackage{enumitem}
\usepackage{bbm} 

\usepackage{xcolor} 

\usepackage{hyperref}
\usepackage{url}
\usepackage{stmaryrd}

\theoremstyle{plain}

\newtheorem{theo}{Theorem}[section] 
\newtheorem{prop}[theo]{Proposition}
\newtheorem{lemme}[theo]{Lemma}
\newtheorem{cor}[theo]{Corollary}
\newtheorem{defin}[theo]{Definition}

\theoremstyle{definition}

\newtheorem{rem}[theo]{Remark}

\newcommand{\cc}{\color{black}}

\newcommand{\cd}{\color{black}}

\newcommand{\R}{{\mathbb{R}}}
\newcommand{\N}{{\mathbb{N}}}
\newcommand{\C}{{\mathbb{C}}}

\newcommand{\be}{{\beta}}
\newcommand{\al}{{\alpha}}
\newcommand{\la}{{\lambda}}

\newcommand{\si}{{\sigma}}
\newcommand{\varsi}{{\varsigma}}

\newcommand{\ga}{{\gamma}}

\newcommand{\om}{{\omega}}
\newcommand{\Om}{\Omega}

\newcommand{\Si}{{\Sigma}}
\newcommand{\ep}{\epsilon}

\newcommand{\Ci}{{\mathcal{C}}^{\infty}} 

\newcommand{\op}{\operatorname}
\newcommand{\con}{\overline}
\newcommand{\bigo}{\mathcal{O}}

\newcommand{\biginf}{\mathcal{O}_{\infty}}

\newcommand{\wt}{\widetilde}
\newcommand{\wh}{\widehat}

\newcommand{\piso}{\Psi_{\op{iso}}}

\newcommand{\dheis}{\mathcal{D}_{\op{Heis}}}
\newcommand{\pheis}{\Psi_{\op{Heis}}}

\newcommand{\ptw}{\Psi_{\op{tsc}}}

\newcommand{\psc}{\Psi_{\op{sc}}}

\newcommand{\dis}{\mathcal{D}_{\op{is}}}
\newcommand{\pis}{\Psi_{\op{is}}}

\newcommand{\tpi}{\tilde{\pi}}

\hypersetup{pdfborder={0000}, colorlinks=true, linkcolor=blue, citecolor=red}

\begin{document}

\title{Resolvents of Bochner Laplacians in the semiclassical limit}
\author{L. Charles \footnote{laurent.charles@imj-prg.fr}}

 \maketitle

\begin{abstract}
  We introduce a new class of pseudodifferential operators, called Heisenberg
 semiclassical pseudodifferential operators, to study the space of sections of
 a power of a line bundle on a compact manifold, in the limit where the
 power is large. This class contains the Bochner Laplacian associated to a
 connection of the line bundle, and when the curvature is nondegenerate, its
 resolvent and some associated spectral projections, including generalised
 Bergman kernels. 
 \end{abstract}

\section{Introduction} 

The spectral analysis of the Bochner Laplacians acting on sections of a line bundle with a large
curvature, has many applications ranging from complex geometry to mathematical
physics: holomorphic Morse inequalities and Bergman kernels \cite{MaMa}, dynamical
systems \cite{FaTs15}, geometric quantization  \cite{BoUr96} or
large magnetic field limit of Schrödinger operators \cite{NicSan}, \cite{Leo}, \cite{Kol} to quote
just a few references. In this paper we introduce
an algebra of pseudodifferential operators, shaped to study the
bottom of the spectra of these Laplacian at small scale.

To understand how the scales matter,  let us state two 
Weyl laws corresponding to two different regimes. Let $(M,g)$ be a closed
Riemannian manifold {\cc of dimension $n$} and $L
\rightarrow M$  a
Hermitian line bundle with a connection $\nabla$. For any
integer $k$, the Bochner Laplacian $\Delta_k = \tfrac{1}{2} \nabla^* \nabla$
acting on sections of $L^k$ is an elliptic differential operator. Hence
$\Delta_k$ with domain $\Ci ( M , L^k)$, is essentially
self-adjoint. Its spectrum is a subset of $\R_{\geqslant 0}$ and consists only
of eigenvalues with finite multiplicity.

When $k$ is large,  the structure of this spectrum depends essentially on the curvature $\frac{1}{i} \om $ of $\nabla$.
Since $\nabla $ is assumed to preserve the Hermitian structure, $\om$ is a
real closed $2$-form of $M$.
For the first regime we will need the twisted symplectic
form $\Om$ of $T^*M$ defined by
\begin{gather} \label{eq:twisted_symplectic_form}
\Om = \sum d\xi_i \wedge dx_i + p^* \om 
\end{gather}
where $p$ is the projection $T^*M \rightarrow M$. For
the second one, we will assume that $\om$ is {\cc non-degenerate, so symplectic,}
and compatible with $g$  in the sense that $\om ( x,y) = g
(jx,y)$ for an (almost) complex structure $j$.
The Weyl laws at the energy scales $k^2$ and $k$ are respectively: 
\begin{enumerate}
  \item For any $ \lambda  >0$, the number $N_k (\lambda )$ of eigenvalues of $k^{-2}
\Delta_k$ smaller than $\lambda$ and counted with multiplicities satisfies
\begin{gather} \label{eq:weyl_law_twisted}
N_k ( \lambda) \sim  \Bigl( \frac{k}{2\pi} \Bigr)^n  \op{vol} \{ \xi \in T^*M, \;
  \tfrac{1}{2} | \xi|_x^2 \leqslant \lambda \} \quad \text{ as $k \rightarrow
    \infty$}  
\end{gather}
with the volume computed with respect to the
  Liouville form $\frac{1}{n!} \Om^n$.
\item if $\om$ and $g$ are compatible, then for any $M > 0$ there exists $C>0$
  such that for any $k$, the spectrum of $k^{-1} \Delta_k$ satisfies
\begin{gather} \label{eq:cluster}
  \op{sp} ( k^{-1} \Delta_k ) \cap ]-\infty, M] \subset ( \tfrac{1}{2} + \N
  ) + C k^{-\frac{1}{2}} ]-1, 1[,  
\end{gather}
and for any $m \in \N$,
\begin{gather} \label{eq:dim_cluster_asymptotic}
  \sharp \Bigl( \op{sp} ( k^{-1} \Delta_k ) \cap  \bigl(
   \tfrac{d}{2} + m  +   ]-\tfrac{1}{4} , \tfrac{1}{4} [ \bigr)  \Bigr) \sim
   \Bigl( \frac{k}{2\pi} \Bigr)^{d} {\textstyle{ { m+d-1 \choose m } }}
   \op{Vol} ( M)
   \end{gather} 
   {\cc  as $k \rightarrow \infty$ with $d = \frac{n}{2}$ and $ \op{Vol} ( M) = \frac{1}{d!} 
\int_M     \om^{d}$.}
 \end{enumerate}

The estimate \eqref{eq:weyl_law_twisted} does not appear in the literature,
but it is actually a small variation of the Weyl law of semiclassical
pseudodifferential operators on a compact manifold with semiclassical
parameter $h = k^{-1}$, as we will explain below.
The cluster structure \eqref{eq:cluster} and the estimate
\eqref{eq:dim_cluster_asymptotic} of the number of eigenvalues in each cluster
have been proved in \cite{FaTs15}. In \cite{oim_part1},
\cite{oim_copain}, it is shown that this eigenvalue number is given by
a Riemann-Roch formula when $k$ is sufficiently large and {\cc generalisations of
\eqref{eq:cluster},  \eqref{eq:dim_cluster_asymptotic} are proved under the
assumption that $\om$ is symplectic but not necessarily compatible with $g$.}

These results rely on a good description of the resolvents $(k^{-\epsilon}  \Delta_k - z)^{-1}$,
which allows to study the $f ( k^{-\epsilon} \Delta_k)$ for $f$ a smooth compactly supported
function, where $\epsilon =2$ or $1$  according to the regime. In the first case,
$f(k^{-2} \Delta_k)$ and $(k^{-2} \Delta_k -z)^{-1}$ are semiclassical twisted pseudodifferential
operators, a class of operators which was introduced in \cite[Chapitre
4]{oim_these}  and used recently in \cite{Colin}.
The local theory of these operators is exactly the same as the one of the standard  semiclassical pseudodifferential
operators with $h = k^{-1}$, whereas the global theory involves the twisted symplectic form
\eqref{eq:twisted_symplectic_form}. Similar operators have been
studied as well on $\R^n$ under the name of magnetic pseudodifferential
operators \cite{MaPu}, cf. \cite{Lein} for an introduction and many
references. 

For the second regime, the theory is much less developed. The proof of \eqref{eq:cluster} is based on an approximation of the resolvent
of $k^{-1} \Delta_k$, which is  obtained by gluing local resolvents of some model
operators deduced from $\Delta_k$ by freezing the coordinates in an appropriate
way. The main difficulty in this construction is that the number of terms we
have to glue increases with $k$. Still these approximations have been used
successfully to prove \eqref{eq:cluster} and later to describe the spectral
projections onto the various clusters as generalizations of Bergman kernels
\cite{oim_part1}, \cite{oim_copain}, \cite{Yuri_1}.

In this paper, we introduce a new class of pseudodifferential operators, and
we prove that it contains the resolvent of $k^{-1} \Delta_k$ when $\om$ is
nondegenerate, as well as the spectral projector
associated to the cluster \eqref{eq:cluster}. This operator class is in a
sense a semiclassical version of the class of Heisenberg pseudodifferential operators, which
has been developed for the study of the $\overline{\partial}_b$-complex
\cite{BeGr}, \cite{EM}. For this reason, we call our operators semiclassical
Heisenberg pseudodifferential operators.

Just as the usual Heisenberg
operators have a symbol of type $(\frac{1}{2}, \frac{1}{2})$, the
semiclassical ones have an exotic symbol of type $\frac{1}{2}$: at
each derivative, a power of $h^{\frac{1}{2}}$ is lost. Recall that the semiclassical operators with a
symbol of type $\frac{1}{2}$ form a limit class: they are closed under
product, the usual operator norm estimates hold, but the standard expansions in
the symbolic calculus do not hold because all the terms in these expansions have the same magnitude, cf. for
instance \cite[Proposition 7.7, Theorem 7.9 and Theorem 7.11]{DiSj}. As we
will see in the sequel \cite{oim_Heis} of this paper, the semiclassical Heisenberg operators
are closed under product. However the composition of their principal symbols, which
are functions on $T^*M$, is not the usual product. Instead it is a
fiberwise product, whose restriction to each fiber $T_x^*M$ depends on the curvature
$\om_x$: when $\om_x$ is nondegenerate, it is essentially the Weyl product
whereas when $\om_x =0$, this is the usual function product. So in general the  principal symbol composition  is not commutative. 

Our definition of the Heisenberg operators is based on the usual semiclassical
pseudodifferential operators, from which we deduce easily many of their
properties, except what regards their composition. In this paper, we only
address the Heisenberg composition of differential operators with
pseudodifferential ones, because this suffices to show that the 
resolvent and cluster spectral projectors of $k^{-1} \Delta_k$ are Heisenberg
operators. We have also included an exposition of the theory of
twisted pseudodifferential operators, because this is not really standard
material, and this helps to understand  the specificity of the Heisenberg
pseudodifferential operators. In the remainder of the introduction, we state
our main result.

\subsection*{Twisted pseudodifferential operators} 

As the usual Laplace-Beltrami operator, $\Delta_k$ has the following local
expression in a coordinate
chart $(U,x_i)$  
$$ k^{-2} \Delta_k =  \tfrac{1}{2 \sqrt g }  \sum_{j, \ell  =1 }^{n}  \pi_j
g^{j \ell}
\sqrt g \, \pi_{\ell} $$
where for any {\cc $j =1, \ldots, n$}, $\pi_j$ is the dynamical moment $ \pi_j :=  (ik )^{-1} \nabla ^{L^k}_j$ with $ \nabla ^{L^k}_j$  the covariant
derivative of $L^k$ with respect to $\partial_{x_j}$.

Let $s \in \Ci ( U, L)$ be such that $|s| =1$ and let $\beta  \in \Om^1 (U,
\R)$ be the corresponding connection 1-form, $\nabla
s = -i \beta \otimes s$.  Identifying $\Ci ( U, L^k) $ with $\Ci (U)$ through
the frame $s^k$,
$$\pi_j = (ik)^{-1} \partial_{x_j} - \beta_j
$$ where
$\beta_j(x)  = \beta(x) ( \partial_{x_j})$. 
Under these local identifications, the $\pi_j$'s and consequently $k^{-2} \Delta_k$ are semiclassical
differential operators with $h =k^{-1}$, their symbols being respectively $\xi_j - \beta_j$ and $\frac{1}{2}
| \xi - \beta (x)|^2$.

These symbols become independent of $s$ if we pull-back
them by the momentum shift $T_\beta : T^*U \rightarrow T^*U$, $(x, \xi) \rightarrow (x ,
\xi + \be (x))$.  The same method will be used to define the symbol of a twisted
pseudodifferential operator, cf. \eqref{eq:def_princ_symbol}. Similarly, the shift $T^*_\be$ is used in classical mechanics
to write the motion equations of a particle in a magnetic field in an invariant way: $T_\be$ sends the
Hamiltonian $\frac{1}{2} |\xi - \beta|^2$ with the standard symplectic form $\sum d\xi_i
\wedge dx_i$ to the usual kinetic energy $\frac{1}{2} |\xi|^2$ with the
twisted symplectic form $\Om$.

Let us briefly introduce the twisted pseudodifferential operators of $L$,
the complete definition will be given in Section \ref{sec:twist-pseud-oper}.
Let us start with the residual  class. An operator family 
\begin{gather} \label{eq:op_family_intro}
P = (P_k :
\Ci ( M, L^k) \rightarrow \Ci (M, L^k), \; k \in \N)
\end{gather}
belongs to $k^{-\infty} \Psi^{-\infty}(L) $ if for each $k$, the Schwartz
kernel of $P_k$ is smooth, its pointwise norm is in $\bigo ( k^{-\infty})$
uniformly on $M$ and the same holds for its successive derivatives.  For
any
local data $\delta =(U, s, \rho)$ consisting  of an open set $U$ of $M$, a
section 
$s \in \Ci ( U, L)$ such that $|s|=1$ and a function $\rho \in \Ci_0(U)$, we
define the local form of a family $P$ as in \eqref{eq:op_family_intro} as:   
\begin{gather} \label{eq:op_local_data}
P^\delta_k  : \Ci( U) \rightarrow \Ci (U), \qquad    (P^\delta_k  f) s^k = \rho P_k ( \rho f
s^k ), \qquad k \in \N
\end{gather}

A twisted pseudodifferential operator of order $m$ is a family $P$ as in \eqref{eq:op_family_intro}
such that for any $\rho_1$, $\rho_2 \in \Ci ( M)$ with disjoint supports,
$(\rho_1 P_k \rho_2)$ belongs to  $k^{-\infty} \Psi^{-\infty} (L)$ and for any local
data $\delta =(U, s, \rho)$,  $P_k^\delta  = Q^\delta_ {k^{-1}} $ where $(Q_h^\delta  :\Ci(U) \rightarrow \Ci (U),
\; h \in (0,1])$ is a semiclassical pseudodifferential operator of order $m$.
So in terms of coordinates $(x_i)$ on $U$  the Schwartz kernel of $P_k^\delta$ has the form 
\begin{gather} \label{eq:kernel_twisted_pseudo}
P_k^\delta (x,y)  = \Bigl( \frac{k}{2 \pi} \Bigr)^n \int e^{ik \xi \cdot (x - y) } a (k^{-1} ,
\tfrac{1}{2} (x+y), \xi) \; d\xi
\end{gather}
for some semiclassical polyhomogeneous symbol $(a (h, \cdot ) \in \Ci ( U
\times \R^n), \; h \in (0,1])$ of order $m$.
The (principal) symbol of $P$ is the function  $\si \in \Ci (T^*M)$ such that for any local
data $\delta$ as above
\begin{gather} \label{eq:def_princ_symbol}
a(h, x, \xi + \beta (x) ) = \rho(x) \si (x, \xi) +
\bigo ( h),
\end{gather}
where $\beta \in \Om ^1 (U, \R)$ is the
connection one-form of $s$.



Examples of twisted (pseudo)differential operators of respective order $0$,
$1$ and $2$ are the multiplications by the functions of $\Ci (M)$,  the covariant derivatives $(ik)^{-1} \nabla^{L^k}_X$ where $X \in \Ci ( M, TM)$, the symbol being 
 $(x, \xi) \rightarrow \langle \xi, X(x)
\rangle$, and the normalised Laplacian  $k^{-2}
\Delta_k$, its symbol is $(x, \xi ) \rightarrow \frac{1}{2} |\xi|_x^2$.

It is not difficult to adapt the standard results \cite{Cdv}, \cite{zworski} on the resolvents of
elliptic operators and their functional calculus in this setting. Denoting by
$\ptw^m(L)$ the space of twisted pseudodifferential operators of order $m$, we have
\begin{itemize} 
 \item for any $z \in \C \setminus \R_{\geqslant 0}$, the resolvent $( z- k^{-2}
   \Delta_k) ^{-1} $ belongs to $\ptw^{-2} ( L) $ and its symbol is  $( z - \frac{1}{2} |
  \xi|^2 )^{-1} $. 
\item for any $ f \in \Ci_0(\R)$, $ f (k^{-2} \Delta_k) $ belongs to
  $\ptw^{-\infty} (L) $ and its symbol is  $f ( \frac{1}{2} |
  \xi|^2 )$. 
\end{itemize}
In particular $\op{tr} f (k^{-2} \Delta_k) = \bigl( \frac{k}{2 \pi} \bigr)^n
\int_{T^*M} f (\frac{1}{2} |
  \xi|^2 ) \frac{1}{n!}  \Om^n  + \bigo (k^{-1+n})$.  
The Weyl law \eqref{eq:weyl_law_twisted} follows. 

\subsection*{Heisenberg pseudodifferential operators} 

A (semiclassical) Heisenberg pseudodifferential operator of $(L, \nabla)$ of
order $m$ is by definition a family $P$ of operators of the form \eqref{eq:op_family_intro}
such that for any $\rho_1$, $\rho_2 \in \Ci(M)$ with disjoint supports,
$(\rho_1 P_k  \rho_2)$ belongs to $k^{-\infty} \Psi^{- \infty} (L)$ and  for any local data $\delta = ( U, s, \rho)$ as above with a
coordinate set $(x_i)$ on
$U$, the Schwartz kernel of $P_k^\delta$ has the form
\begin{gather} \label{eq:intro_Def_Heisenberg}
  e^{ i k \be \bigl( \tfrac{x+y}{2} \bigr) \cdot (x-y) } \Bigl(
  \frac{\sqrt{k}}{ 2 \pi } \Bigr)^n \int_{\R^n} e^{i \sqrt{k} \; \xi \cdot
    (x-y) } a ( k^{-\frac{1}{2}} ,\tfrac{1}{2} ( x+y), \xi ) \; d\xi 
\end{gather}
where
  \begin{enumerate}
    \item[-]   $\be = \sum \beta_i(x) dx_i $ is the connection one-form of
      $s$ defined as above and   
 viewed as the  $\R^n$-valued function  $x \rightarrow (\be_1(x), \ldots ,
 \be_n(x))$
 \item[-] $ (a ( h , \cdot ) \in \Ci ( U  \times \R^n), \; h \in (0,1])$ is a
   semiclassical polyhomogeneous symbols of order $m$,  so in particular $\partial_x^{\al}
   \partial_\xi^\beta a = \bigo_{\al, \be} ( \langle \xi \rangle ^{m- | \beta|})$ and  $a
   \sim \sum_{\ell =0 } ^\infty h ^\ell a_\ell ( x, \xi) $
 with polyhomogeneous coefficients $a_{\ell}$ of order $m- \ell$.
\end{enumerate}
As we see, the Schwartz kernel \eqref{eq:intro_Def_Heisenberg} is the product
of an oscillatory factor depending on the frame $s$ with the Schwartz kernel
of a semiclassical operator where the semiclassical parameter is
$k^{-\frac{1}{2}}$. We will prove that this formula is consistent with change
of frame and that we can define a (principal) symbol $\si \in \Ci ( T^*M)$
such that 
for any local data as above,
$$a( h, x, \xi ) = \rho (x) \si (x, \xi) + \bigo
(h). $$
Let us denote by $\pheis^m ( L, \nabla)$ the space of Heisenberg
pseudodifferential operators of order $m$, and by $\pheis^{-\infty} (L, \nabla)$ the
intersection  $\cap_m \pheis^m ( L, \nabla)$.

To state our main result, we need to introduce some symbols  $R_{d,z}$ and
$\pi_{d,E}$. Recall that 
 for any tempered distribution $ a \in \mathcal{S}' ( \R_s^d \times \R_\varsigma^d)$, the Weyl
quantization of $a$ is the operator $ a^w  :
\mathcal{S} ( \R^d) \rightarrow \mathcal{S}'( \R^d)$ with Schwartz kernel at
$(s,t) \in \R^d \times \R^d$: $$ (2
\pi)^{-d} \int e^{i \varsigma \cdot (s -t ) } a( \tfrac{1}{2} (s+t) , \varsigma ) \; d
\varsigma.$$
{\cd{Here the coordinate $s \in \R^d$ should not be confused with the frame $s$
introduced previously.}} 
The (quantum) harmonic oscillator is $H^w$ with $H (s, \varsigma) = \tfrac{1}{2} \sum_{i=1}^d
(s_i^2 + \varsigma_i^2)$. As an operator of $L^2 ( \R^d )$ with domain the Schwartz
space, $H^w$ is essentially self-adjoint with spectrum $\op{sp} H^w = \frac{d}{2} + \N$.
Then for any $ z \in \C \setminus \op{sp} H^w$ and $E \in \op{sp} H^w$,
$R_{d,z}$ and $\pi_{d,E}$ are the tempered distributions of $\R^{2d}$
such that 
\begin{gather} \label{eq:symb_non_com}
  ( H^w - z )^{-1} = R_{d,z} ^w , \qquad 
1_{\{ E \} } ( H^w) = \pi^w_{d,E}
\end{gather} 
By  Weyl calculus, $R_{d,z}$ belongs to the symbol class
$S^{-2} ( \R^{2d})$ and $R_{d,z} = (H - z )^{-1}$ modulo $S^{-3} (\R^{2d})$.
{\cc The definition of the symbol class $S^m $ will be recalled in Section
\ref{sec:twist-pseud-oper}.} Moreover $\pi_{d,E} $ belongs to the Schwartz space,
being the Weyl symbol of an orthogonal projector onto a finite dimensional subspace of $\mathcal{S} (
\R^d)$. The analytic Fredholm theory can be developed in this setting and it
says that the function $z \rightarrow R_{d,z} $ with values in $\Ci ( \R^{2d})$, or better
the symbol space $S^{-2} ( \R^{2d})$, is meromorphic on $\C$ with simple poles at
$ \frac{d}{2} + \N$ whose residues are the $\pi_{d,E}$.

{\cd{In Theorem \ref{theo:main_intro}, we will view $H$, $R_{d,z}$ and $\pi_{d,E}$ as functions on
$T^*_xM$, for any $x \in M$ as follows: we assume that $\om$ and $g$ are compatible so
that $ n= 2d$ with $d \in \N$.  We choose an
orthosymplectic basis $(e_i,f_i)_{i=1}^d$ of $T_xM$, that is $(e_i,f_i)$ is an orthonormal
basis and for any $i,j$,
$$\om (x) ( e_i , e_j) =0 =\om (x) (f_i, f_j), \qquad \om (x) ( e_i, f_j) =
\delta_{ij} .$$ Let $(s_i,\varsi_i)$ be the associated
coordinates of $T^*_xM$, so   $s_i ( \xi) := \xi (e_i)$ and $\varsi_i ( \xi )
:= \xi (f_i)$ for any $\xi \in T_x^*M$. 
Then any function $f:\R^{2d} \rightarrow \C $  identifies with  the
function of $T_x^*M$
\begin{gather} \label{eq:identif}
\xi \in T^*_x M \rightarrow f( s( \xi), \varsigma ( \xi)).
\end{gather}
In
particular $H( s ( \xi ), \varsigma ( \xi)) = \frac{1}{2} |\xi|^2$ because
$(e_i,f_i)$ is orthonormal. The fact that $ R_{d,z} (s(\xi),
\varsigma ( \xi))$ and $\pi_{d,E} ( s( \xi), \varsigma ( \xi))$ are independent of
the choice of the basis $(e_i,f_i)$ follows from the symplectic invariance of
the Weyl quantization.

It follows as well from the symplectic invariance of Weyl quantization and the $O(n)$ invariance of
$H$ that $R_{d,z}$ and $\pi_{d,E}$ are radial functions. A computation from Mehler formula leads to \cite{De}
\begin{gather} \label{eq:Derez}
R_{d,z} = \int_0^1 ( 1 - \tfrac{1}{2}s ) ^{\frac{d}{2} - z - 1 } ( 1 +
\tfrac{1}{2}s ) ^{ \frac{d}{2} + z -1 } e^{- s H }  ds , \quad  \text{ if }
\op{Re} z < d
\end{gather}
We can also compute $\pi_{d,E}$ in terms of Laguerre polynomials
\cite{Un} :
\begin{gather} \label{eq:Un}
\pi_{d,E} (\xi) = 2^d (-1) ^m
e^{-|\xi|^2}L_m^{d-1} ( 2 | \xi |^2), \qquad \text{ where } m = E -
\tfrac{d}{2}
\end{gather}
and $L_m ^\al (x) = \frac{1}{m!} e^x
x^{-\al} \partial_x^m ( e^{-x} x^{m+ \al}).  
$}}

\begin{theo} \label{theo:main_intro}
  Assume $\om$ and $g$ are compatible so that $n = 2d$ with $d \in \N$. Then
  \begin{enumerate}
    \item 
For any $ z \in \C \setminus ( \tfrac{d}{2} + \N)$, there exists $Q(z) \in
\pheis^{-2} ( L, \nabla)$ such that
\begin{enumerate}
  \item[-] $(k^{-1} \Delta_k - z ) Q_k(z) \equiv \op{id}$ and $  Q_k(z)
    (k^{-1} \Delta_k - z ) \equiv \op{id}$ mod $k^{-\infty} \Psi^{-\infty}
    ( L)$. 
\item[-] $(k^{-1} \Delta_k - z ) Q_k(z) = Q_k(z)
    (k^{-1} \Delta_k - z ) = \op{id}$ when $k$ is large. 
  \item[-] the {\cc principal} symbol of $Q (z)$ restricted to $T_x^*M$ is the Weyl symbol
    $R_{d,z}$ of
    the resolvent of the harmonic oscillator with symbol $\tfrac{1}{2} |\xi|^2$, cf. \eqref{eq:symb_non_com}.  
\end{enumerate} 
\item 
  For any $E \in \frac{d}{2} + \N$, the spectral projector family
  $$(1_{[E-
  1/2,E + 1/2]} ( k^{-1} \Delta_k),\; k \in \N)$$ belongs to $\pheis^{-\infty} (
L, \nabla)$. The restriction of its {\cc principal} symbol to $T_x^*M$ is the Weyl symbol
$\pi_{d,E}$ of
the spectral projector on the $E$-eigenspace of the harmonic oscillator with symbol $\tfrac{1}{2} |\xi|^2$, cf. \eqref{eq:symb_non_com}
\end{enumerate}
\end{theo}

{\cc The Weyl estimates \eqref{eq:dim_cluster_asymptotic} are a consequence of Theorem
\ref{theo:main_intro}. Indeed, when $\om$ is symplectic, we have the following
equivalent for the trace of an Heisenberg pseudodifferential
operator $P$ of $(L, \nabla)$ of order $-\infty$ with symbol $\si$
$$ \op{tr} P_k = \Bigl( \frac{k}{2 \pi} \Bigr)^d \int_M \op{tr} \si^w (x) \;d
\mu_M (x) \qquad \text{with } \mu_M = \tfrac{1}{d!} \om^d .$$
For $P = 1_{[E-  1/2,E + 1/2]} ( k^{-1} \Delta_k)$, $\op{tr} P_k$ is the
dimension of the cluster and $\op{tr} \si^w (x)$ is the rank of $\pi_{d,E}$,
that is $ { m+d-1 \choose m }$ for $ E = \frac{d}{2} +m$.}

  The proof of the first part of Theorem \ref{theo:main_intro} is an adaptation
of the standard
parametrix construction of an elliptic pseudodifferential operator, the main
change being in the symbolic calculus:  if $(P_k)$ belongs to $\pheis^m (L, \nabla)$ and has symbol $\si$, then $ ( k^{-1}
\Delta_k P_k)$ belongs to {\cc $\pheis^{m
    +2} ( L , \nabla)$} and its symbol 
restricted to $T_x^*M$ is the Weyl product of $\tfrac{1}{2} |\xi|_x^2$ and $ \si
(x, \cdot)$. By Weyl product, we mean the product of symbols in Weyl
quantization,  and the identification of functions of $T_x^*M$ with symbols is
done through \eqref{eq:identif}. This explains how  the symbols $R_{d,z}$
and $\pi_{d,E}$ appear.

 A remarkable fact is that the proof of the second assertion of Theorem
\ref{theo:main_intro} is a direct application of Cauchy formula for the
spectral projector of an operator with compact resolvent. This part is much
simpler than the proof that $f( k^{-2} \Delta_k) \in \ptw^{-\infty} ( L)$,
even with the modern approach through Helffer-Sjöstrand formula.

\begin{rem} $ $

{\cc 1. In the sequel \cite{oim_Heis}, we will prove that the Heisenberg
pseudodifferential operators form an algebra. The proof is rather long and
technical, which explains why we choose to
postpone it to another paper. In the current paper, we will merely consider
compositions of differential Heisenberg operators with pseudodifferential
ones. The corresponding symbolic calculus is already nontrivial, and it is sufficient to construct the resolvent of $k^{-1} \Delta_k$.}

2. The Schwartz kernel of the cluster spectral projectors was
described in \cite{Yuri_1} and \cite{oim_copain} as a generalisation of the Bergman
kernel. {\cc The paper \cite{oim_copain} was based on a particular algebra  $\mathcal{L} (
\C)$, introduced in \cite{oim_part1}, containing the cluster spectral projectors and the
associated Toeplitz operators. This algebra  $\mathcal{L}(\C)$ was used to compute  the dimension of each cluster and to develop the theory
of Toeplitz operators. The algebra $\pheis^{-\infty} (L, \nabla)$ could be
used equivalently once the composition theorem is proved.

The relation between
both algebras is simply that
$\mathcal{L}(\C)$ is a subalgebra of $\pheis^{-\infty} (L, \nabla)$. Each 
algebra has its own advantages. The Schwartz kernels of the operators of
$\mathcal{L}( \C)$ are described without oscillatory integrals; and their
asymptotics, given directly in their definition, is more precise than what can be
obtained for the Heisenberg pseudodifferential operators. 
The advantage of considering  Heisenberg
 operators is merely that $\pheis (L, \nabla)$ contains the Laplacian and its
resolvent, so that the cluster spectral projectors and the resolvent are described
uniformly. The proof in \cite{oim_copain} that the cluster spectral projector
belongs to $\mathcal{L} ( \C)$ used an approximation of the resolvent, which
was not satisfactory as explained earlier. }
\qed \end{rem}

\subsubsection*{Outline of the paper} 
In Section \ref{sec:twist-pseud-oper}, we introduce notations and basic analytical tools to
address the large $k$ limit of the space of sections of the $k$-th power of $L$, including  the theory of semiclassical twisted
pseudodifferential operators with their Sobolev spaces. The study of 
Heisenberg pseudodifferential operators starts in Section
\ref{sec:heis-semicl-oper} {\cc with the global intrinsic definition of their
Schwartz kernels. We then consider local expression of these kernels in
Section \ref{sec:local-expressions} and establish mapping properties in Section \ref{sec:mapping-properties}.} In Section \ref{sec:isotropic-algebra}, we introduce the symbol
product, which is then used in Section \ref{sec:heis-diff-oper} for the composition of differential
operators with pseudodifferential operators. This is applied to resolvents
and spectral projections in Section \ref{sec:resolvent}. In Section \ref{sec:auxiliary-bundles}, we explain how we can
add auxiliary bundles to the theory, which provides some important examples.

\subsubsection*{Acknowledgment}
I would like to thank the anonymous referee for a careful reading of the paper and useful
comments. I would like to thank as well Clotilde Fermanian Kammerer, Colin Guillarmou and
Thibault Lefeuvre for
useful discussions.

\section{Twisted pseudodifferential operators} \label{sec:twist-pseud-oper}

\subsubsection*{Symbols} 
We will use 
the class of semiclassical polyhomogeneous symbols
 introduced in \cite[Section E.1.2]{DiZw}, cf. also \cite[Section 6.1]{Cdv}.
Let $V$ be an open set of $\R^p$ and $m \in \R$. For any $\xi \in \R^n$, let
$|\xi|$ and $\langle \xi \rangle$ be the Euclidean norm and Japanese
bracket, so $| \xi |^2 = \sum \xi_i ^2$, $\langle \xi \rangle ^2 = 1 + |
\xi|^2$. Let $S^m ( V,
\R^n)$, $S^m_{\op{ph}} (V, \R^n)$ and $S^m_{\op{sc}} (V, \R^n)$ be the spaces of
symbols (resp. polyhomogeneous symbols, semiclassical polyhomogeneous
symbols)  of order $m$. By definition
\begin{itemize} 
\item $S^m (V, \R^n)$ consists of the families $(a(h, \cdot),
\; h \in (0,1])$ of $\Ci ( V \times \R^n)$ such that for any compact set $K$
of $V$, $\al\in \N^p , \be \in \N^n$, there exists $C>0$ such that 
$$ |\partial_x^\al \partial_\xi ^\be a (h,x,\xi) | \leqslant C \langle \xi
\rangle^{m- | \be|}, \qquad \forall x \in K, \; \xi \in \R^n, \; h \in
(0,1] $$
\item $b \in S_{\op{ph}}^m (V, \R^n)$ if $b \in S^m (V, \R^n)$, $b$  is independent
  on $h$ and for
every $N$,  $b = \sum_{j=0}^{N-1} b_j \mod S^{m-N} (V, \R^n)$ with coefficients $b_j \in \Ci ( V \times \R^n)$ such that  $b_{j} (x, t\xi) = t^{m-j} b_{j} (x,\xi)$ when
$|\xi| \geqslant 1$ and $t\geqslant 1$.
\item $a \in S^m_{\op{sc}} (V, \R^n)$ if $a \in S^m (V, \R^n)$ and {\cc for every
  $N$,} $ a = \sum_{\ell=0}^{N-1} h^{\ell} a_{\ell} $ mod $
 h ^{N} S^{m-N} (V, \R^n)$ for some coefficients $a_{\ell} \in S_{\op{ph}}^{m-\ell} ( V ,
\R^n)$. 
\end{itemize}
More generally these definitions make sense for a real vector bundle $E \rightarrow
N$ instead of the product $V \times \R^n$.  We denote by $S_*^m (N, E)$ the
corresponding spaces and set
\begin{gather} \label{eq:symb_fibre}
S^{\infty}_* (N, E) = \bigcup_m S^m_* (N,E), \qquad  S^{-\infty}_* (N, E) =
\bigcap_m S^m_* (N,E) 
\end{gather}
for
$* = \emptyset, \op{ph}, \op{sc}$.
An easy remark is that for any section $u$
of $E$, the translation $T_u $ of $\Ci ( E, \C)$ given by
\begin{gather} \label{eq:translation}
T_uf ( x,v ) = f(
x, v - u(x))
\end{gather}
preserves $S_*^m ( N, E)$.  When $V$ is reduced to a point, we
set $S^m_* ( \R^n) := S^m_* ( \{ \cdot \}, \R^n)$. Moreover,
$S^{-\infty}_{\op{ph}} ( \R^n) $ is the Schwartz space $\mathcal{S} ( \R^n)$. 

\subsubsection*{Negligible families} 

We say that a family $(f_h , h \in (0,1])$ of $\Ci(N)$ is negligible, and
we write $f_h = \bigo_{\infty} ( h^{\infty})$, if all its $\Ci$-seminorms are
in $\bigo ( h^{\infty})$. This definition is meaningful if $(f_h)$ is only
defined for $ h \in D$ where $D$ is any subset of $(0,1]$ whose closure
contains $0$. 

Let $L \rightarrow M$ be a Hermitian line bundle and $A \rightarrow M$ a
complex vector bundle with rank $r$. {\cc Introduce an open set $U$ of $M$ with a
frame $(a_i \in \Ci  (U, A),\; i=1, \ldots , r )$ of $A$ and a frame $s \in \Ci
(U, L)$ with pointwise norm
$|s|=1$. Since $L$ has rank one, a frame of $L $ is merely a local section
which does not vanish anywhere. The $k$-th power $s^k$ is a frame of $L^k$. 
Any smooth section of $L^k \otimes A$ over $U$
has the
form $\sum_{i=1}^r f_i  s^k \otimes a_i $ with coefficients $f_ i \in \Ci (U)$. 

A family $ (t_k \in \Ci ( M , L^k \otimes A), \; k \in \N)$ is said to be {\em negligible} if
for any choice of $U$, $(a_i)_{i=1}^r$ and $s$ as above, we have 
\begin{gather} \label{eq:loc}
  t_k = \sum_{i=1}^r f_{i, k^{-1}} s^k \otimes a_i \text{  on } U \qquad
  \text{ with }  
  f_{i,h} \in \bigo _{\infty} (
h^{\infty}). 
\end{gather} 
We denote by $\biginf (k^{-\infty})$ the space of negligible
families.

Notice that the condition \eqref{eq:loc} is independent of the choice of the
frames:  let $(a_i)_{i=1}^r$ and $(\wt{a}_i)_{i=1}^r$ be two frames of $A$
over $U$ and $s, \wt{s}
\in \Ci (U,L)$ two sections such that $|s| =|\wt{s}|=1$. Then $\wt{a}_j = \sum_{\ell=1}^r
g_{j\ell} a_\ell$ with smooth coefficients $g_{ij}$ and restricting $U$ if necessary, $\wt{s} = e^{i \varphi } s$
with $\varphi \in \Ci (U, \R)$. Write  $t_k =  \sum_\ell
f_{\ell, k^{-1}} s^k \otimes a_\ell =  \sum_j
\wt{f}_{j, k^{-1}} \wt{s}^k \otimes \wt{a}_j$. Then 
\begin{gather} 
f_{\ell, h} = e^{ih^{-1}  \varphi } \sum_{j}  g_{j\ell} \wt{f}_{j,h} .
\end{gather}
and we see that  $\wt{f}_{j,h} = \bigo_\infty ( h^{\infty} )$ for any
$j$ implies that $f_{\ell, h} = \bigo_\infty ( h^{\infty} )$ for any $\ell$.   
}

Let $P$ be a family of operators 
\begin{gather} \label{eq:familyp}
 P = (P_k : \Ci ( M , L^k) \rightarrow \Ci ( M , L^k), \; k \in \N
),
\end{gather}
The Schwartz kernel of each $P_k$ is a section of $(L^k \boxtimes \con{L}^k)
\otimes ( \C_M \boxtimes |\Lambda | ( M))$, where we denote by $\boxtimes$ the
external tensor product of vector bundles, by $\C_M$ the trivial line bundle
over $M$ and by $ |\Lambda | ( M)$ the density bundle.    Since $L^k \boxtimes \con{L}^k = ( L \boxtimes
\con{L})^k$, the previous definition of a negligible family applies to the
family $(P_k)$ of Schwartz kernels.

We denote by $k^{-\infty} \Psi^{-\infty}
(L)$ the space consisting of operator
families of the form \eqref{eq:familyp} such that each $P_k$ is smoothing with a
Schwartz kernel family  in $\biginf ( k^{-\infty})$.  As we will see,
$k^{-\infty} \Psi^{-\infty} (L)$ is both the residual space of twisted
pseudodifferential operators and of Heisenberg pseudodifferential operators.

\subsubsection*{Semiclassical pseudodifferential operators}

Let $\Psi^m_{\op{sc}} (M)$ be the space of semiclassical
pseudodifferential operators of order $m$ acting on smooth functions of $M$.
By definition $P \in \Psi^m_{\op{sc}} (M)$ is a family of operators $(P_h: \Ci (M)
\rightarrow \Ci (M), h \in (0,1])$ with a Schwartz kernel $K_h (x,y)$
satisfying for any $\rho \in \Ci (M^2)$, 
\begin{enumerate} 
\item if $\op{supp} \rho \cap \op{diag} M = \emptyset$, then $\rho K_h$
  is smooth and negligible.   
\item if $\op{supp} \rho \subset U^2$ where $(U,x_i)$ is a coordinate chart
  of $M$, then on $U^2$ 
\begin{gather} \label{eq:def_pseudo_semiclassique}
  (\rho K_h) (x,y) = ( 2 \pi h )^{-n} \int e^{i h^{-1} \xi \cdot (x-y) } a( h, x,y , \xi )  \;
d\xi 
\end{gather}
with $a \in S_{\op{sc}}^m (U^2, \R^n)$.
\end{enumerate}
Here and in the sequel, when the Schwartz kernel is written in a coordinate
chart, we implicitly use the density $|dx_1 \ldots dx_n|$. 
The principal symbol of $P$ is the function $\si \in S_{\op{ph}}^m ( M, T^*M)$
such that $ a ( h, x,x, \xi) = \rho (x) \si ( x, \xi) +
\bigo ( h)$ on $U$.  


\subsubsection*{Twisted pseudodifferential operators} 
Let $L \rightarrow M$ be a Hermitian line bundle.


\begin{defin} \label{def:twisted_pseudo} 
  A semiclassical twisted pseudodifferential operator $P$ of $L$
  is a family having the form \eqref{eq:familyp} such that for any $\rho \in \Ci(M^2)$,
  
  \begin{enumerate}
    \item if $\op{supp} \rho \cap \op{diag} M = \emptyset$, then $\rho
      P_k$ is smooth and negligible. 
\item  if $\op{supp} \rho \subset U^2$ where $(U,x_i)$ is a coordinate chart
  of $M$ and $s \in \Ci ( U, L)$ is such that $|s| =1$, then on $U^2$ 
\begin{gather} \label{eq:local_twisted_schwartz}
  ( \rho P_k) (x,y) = \Bigl( \frac{k}{2 \pi} \Bigr)^{n} \int e^{i k \xi \cdot (x-y) } a( k^{-1}, x,y , \xi )  \;
d\xi  \; s^k(x) \otimes \con{s}^k(y) 
\end{gather}
with $a \in S_{\op{sc}}^m (U^2, \R^n)$.
    \end{enumerate}
  \end{defin}
To understand the dependence of the oscillatory integral with respect to the
choice of the frame $s$, consider a new frame $t = e^{i \varphi} s$ where $\varphi \in \Ci(U
, \R)$. Then $\varphi(x) - \varphi (y) = \sum_j \psi_j (x,y) (x_j - y_j)$ where
$\psi_j \in \Ci ( U^2, \R)$ is such that $\psi_j (x,x) = \partial_{x_j} \varphi (
x)$. Using that  $s(x) \otimes \con{s} (y) = \exp ( - i \psi (x,y) \cdot (x-y))
t(x) \otimes \con{t} ( y) )$ and changing the variable $\xi$ into $ \xi  +
\psi (x,y)$, we obtain
\begin{xalignat}{2}  \label{eq:thecalcul}
  \begin{split} 
 & \int e^{i k \xi \cdot (x-y) } a ( k^{-1}, x,y , \xi )  \;
d\xi  \; s^k(x) \otimes \con{s}^k(y) \\ 
& = \int e^{i k \xi \cdot (x-y) } a ( k^{-1},
x,y , \xi + \psi (x,y) )  \;
d\xi  \; t^k(x) \otimes \con{t}^k(y)
\end{split}
\end{xalignat}
So multiplying $s$ by $e^{i \varphi}$ amounts to change the amplitude $a$ to $b$
such that $b( h, x, y, \xi) = a ( h , x, y, \xi + \psi (x,y))$. {\cc As noticed in
\eqref{eq:translation}, if $a \in S_{\op{sc}}^m (U^2, \R^n) $, then the same
holds for $b$. }
Moreover, we have  on the diagonal $$b(h,x,x,\xi) = a ( h,x,x, \xi + d\varphi
(x)).$$
Let $\nabla$ be a connection on $L$ preserving the metric. Then  $\nabla s = \frac{1}{i} \beta_s
\otimes s$, where $\beta_s$ is a real one-form of $U$. Observe that $\nabla t = \frac{1}{i} \beta_t \otimes t$ where $\beta_t =
\beta_s - d \varphi$. So we can define the principal symbol as follows. 
\begin{defin}  The principal symbol $\si_{\nabla} (P)$ of $P \in
  \ptw^m (L)$ is the element of $S^m_{\op{ph}} (M, T^*M)$ 
such that for any local data $(\rho, U, s, a)$ as in Definition
\ref{def:twisted_pseudo}, we have
$$a ( h ,x, x , \xi + \beta_s (x)  ) = \rho  (x) \si_\nabla(P) ( x, \xi )  +
\bigo (h) $$
\end{defin}
If $\nabla'$ is another connection of $L$ preserving the metric, then $\nabla'
= \nabla + \frac{1}{i} \al$ with $\al \in \Om^1 ( M , \R)$ and
$\si_{\nabla'} ( P ) (x, \xi)  = \si_{\nabla} ( P ) (x , \xi + \al (x))$.
It is easy to extend the basic properties of pseudodifferential operators to
our setting:  
\begin{itemize} 
\item[-]  If $P\in \ptw^m (L)$, then $\si_{\nabla} (P) =0 $ if and only if $P \in  k^{-1} \ptw^{m-1}  
(L)$. 
\item[-] $\bigcap_m  k^{-m} \ptw^{-m} (L) = k^{-\infty} \Psi^{-\infty} (L)$.
\item[-]   if $P \in \ptw^m(L) $ and $Q \in
  \ptw^p(L)$, then
  \begin{enumerate}
    \item[i)] $(P_kQ_k)$  belongs to $\ptw^{m+ p}(L)$ and its principal symbol is the
      product of the principal symbols of $P$ and $Q$.
      \item[ii)] $ik [ P_k, Q_k]$
belongs to $\ptw^{m+p-1} (L) $ and its {\cc principal symbol is the Poisson
  bracket of the principal symbols of $P$ and $Q$}, with the Poisson structure
dual to the twisted symplectic form \eqref{eq:twisted_symplectic_form} where
$\frac{1}{i} \om$ is the curvature of $\nabla$.
\end{enumerate}
\end{itemize}

{\cc It is possible to define the twisted pseudodifferential operators without
using local frames of $L$, but a single local frame of $L \boxtimes \con{L}$
defined on a neighborhood of the diagonal. This leads to a direct description
of the twisted pseudodifferential operators in terms of the usual
semiclassical pseudodifferential operators. Since we will use later a similar
presentation for Heisenberg pseudodifferential operators, let us explains how
this construction works. 

 Introduce an open neighborhood
$V$ of the diagonal of $M^2$ and a section $F \in \Ci
( V, L \boxtimes \con {L} )$ such that $|F|= 1$ on $V$ and  $F(x,x) = 1 $ for
any $x \in M$, in the sense that $F(x,x) = u \otimes \con{u}$ for any $u \in
L_x$ with norm $1$. We will prove in Lemma \ref{lem:existenceVF} that
such a section exists. Let $\phi \in \Ci_0( V)$ be equal to $1$ on a neighborhood of the diagonal.

\begin{prop}
A family $P$ of the form \eqref{eq:familyp} belongs to $ \ptw^m (L)$ if and only its Schwartz kernel has
the form
\begin{gather} \label{eq:def_alt_twisted_pseudo}
F^k(x,y) \phi (x,y) K_{k^{-1} } (x,y) + \biginf ( k^{-\infty}) 
\end{gather}
where $(K_h , \; h \in (0,1])$ is the Schwartz kernel family of a semiclassical
pseudodifferential operator $Q \in \psc^m ( M)$. If furthermore $\nabla$ is a
connection of $L$ such that the corresponding covariant
derivative of $F$ is zero on the
diagonal, then $\sigma _{\nabla} ( P) = \si (Q)$.
\end{prop}
\begin{proof}
Assume that \eqref{eq:def_alt_twisted_pseudo} holds. Consider a contractible open set $U$
of $M$ with coordinates $(x_i)$ and $s\in \Ci (U, L)$ such that $|s|=1$. We
have
\begin{gather} \label{eq:Fs}
F(x,y) = \exp ( i \varphi (x,y) )
s (x) \otimes \con{s} ( y)
\end{gather}
with $\varphi \in \Ci (U^2, \R)$ vanishing on the
diagonal. So there exists $\psi_j \in \Ci (U^2, \R)$, $j=1, \ldots,n$ such that $\varphi (x,y)  = \sum_j \psi_j (x,y) (x_j - y_j)$. 
Then by the same
computation as in  \eqref{eq:thecalcul}, we deduce the local expression
\eqref{eq:local_twisted_schwartz} of the Schwartz kernel of $P$ from the
local expression \eqref{eq:def_pseudo_semiclassique} of $K$. If moreover the
symbols in \eqref{eq:def_pseudo_semiclassique} and
\eqref{eq:local_twisted_schwartz} are denoted respectively by $a$ and $b$,
then
\begin{gather} \label{eq:b_et_a}
b ( k^{-1}, x, y, \xi ) = a ( k^{-1} , x, y, \xi - \psi (x,y)).
\end{gather}
This
proves that $P \in  \ptw^m (L)$. The converse is similar.

To prove the second claim, assume that $\nabla s = \frac{1}{i} \beta \otimes
s$ and that \eqref{eq:Fs} holds. Then the covariant
derivative of $F$ is
$$ \nabla F = i ( d \varphi - \pi_{\ell}^* \be + \pi_r^* \be ) \otimes F$$
where $\pi_{\ell} $ and $\pi_r$ are the projections $U^2 \rightarrow U$ on the
left and right factor respectively. So the condition that $\nabla F$ is zero
on the diagonal is equivalent to $\beta_x = \sum  \psi _j (x,x) dx_j$, which we
write as before as $\psi (x,x)   = \beta_x$. So by \eqref{eq:b_et_a}, we have  
that $a ( h,x,x, \xi ) = b (h,x,x, \xi + \beta_x)$, which amounts to say that
$\si ( Q) = \si_{\nabla} (P)$. 
\end{proof}} 


\subsubsection*{Semiclassical Sobolev norms}
Let $m \in \R$. Denote by $H^m  (M , L^k)$ the Sobolev space of
sections of $L^k$ of order $m$. Let us give three equivalent definitions of the semiclassical Sobolev norms
of a section $u$ of $L^k$. First the norm of $H^0 (M , L^k) = L^2 ( M, L^k)$
is defined by
\begin{gather} \label{eq:l2norm}
\| u \|_{L^2 ( M , L^k)}^2 = \int_M | u(x)|^2 d\mu (x) 
\end{gather}
where $\mu$ is a volume element of $M$ independent of $k$.
\begin{enumerate} 
\item only for integral exponent $m \in \N$: choose a connection
  $\nabla$ of $L$, vector fields $(X_i)_{i=1}^N$
of $M$ which generates $T_xM$ at each $x$, and set 
$$\| u \|_m := \sum_{| \al | \leqslant m} k^{-|\al|} 
\| \nabla_X^\al u \|_{L^2 ( M , L^k)}$$
where  for any $\al \in \N^N$,
$\nabla_X^\al = \nabla_{X_1}^{\al(1)} \ldots \nabla_{X_N}^{\al ( N)}$.
\item based on local semi-norms: 
for any chart $(U,\chi )$
of $M$, frame $s \in \Ci ( U, L)$ such that $|s|=1$ and $\rho \in
\Ci_o(U)$ we set  
$$ \| u \|_{m,U,\chi,s, \rho} = \| \langle k^{-1} \xi \rangle^m \hat{v} (\xi) 
\|_{L^2 ( \R^n)} \qquad \text{where }  \rho u = (\chi^* v) s^k $$ 
and $\hat{v}$ is the Fourier transform of $v$. Choose a finite family $(U_i, \chi_i, s_i, \rho_i)$ of local data such
that $M$ is covered by the $\{ \rho_i = 1 \}$ and set
$ \| u \|_m := \sum_i  \| u \|_{m,U_i,\chi _i ,s_i, \rho_i}$.  
\item based on twisted pseudodifferential operators: choose $ E \in \ptw^m ( L)$ which is elliptic and invertible for any $k$, and set
$ \| u \| _m = \| Eu \|_{L^2 ( M, L^k)}$. 
\end{enumerate}
The ellipticity condition is as usual that the principal symbol satisfies for
some $C>0$, $|
\si_{\nabla} ( L) (x, \xi )| \geqslant C^{-1} | \xi|^m$ when $| \xi |
\geqslant C$. It does not depend on the choice of $\nabla$.   

We claim that all these norms are equivalent with constants uniform in $k$.
Furthermore for any twisted pseudodifferential operator $P \in \ptw^p (L)$ and
any $m \in \R$, there exists $C$ such that for any $k$,
\begin{gather} \label{eq:soblov_semi_classical}
\| P_k u \|_m \leqslant C \| u \|_{m+p} , \qquad \forall \; u \in \Ci ( M ,
L^k) .
\end{gather}

\section{Heisenberg semiclassical operator}  \label{sec:heis-semicl-oper}

In the introduction, we defined the Heisenberg
pseudodifferential operators by expressing locally their Schwartz kernels as 
oscillatory integrals. Here we will start with a global definition
which has the advantage that we can deduce some basic properties of these
operators directly from the ones of the semiclassical pseudodifferential
operators.  

Let $L \rightarrow M$ be a Hermitian line bundle with a connection $\nabla$
preserving the metric. The line bundle $L \boxtimes \con{L}$  inherits from
$L$ a Hermitian metric and a connection. Its restriction to the diagonal is
the flat trivial bundle with  a natural trivialisation
obtained by sending $ u \otimes \con{v} \in L_x \otimes
\con{L}_x$ to the scalar product of $u$ and $v$. In the sequel we will use a
particular extension of this trivialisation.

\begin{lemme} \label{lem:existenceVF}
  There exist a tubular neighborhood $V$ of the diagonal of $M^2$ and $F \in \Ci ( V, L \boxtimes \con{L})$ such that $|F | = 1$ on $V$ and
  $$ F(x,x) = 1, \quad \nabla F (x,x) = 0 , \quad \nabla_Y \nabla_Y F (x,x) =0
  \qquad \forall x \in M$$
  for any vector field $Y$ of $M^2$ having the form $Y(x,y) = (X(x) , -X(y))$
  with $X \in \Ci ( M , TM)$.
If $ (V',F')$ satisfies the same conditions, then $F =
  F' \exp ( i \psi )$ where $\psi \in \Ci (V \cap V', \R)$ vanishes to third order along the diagonal 
\end{lemme}

\begin{proof} Consider  more generally a closed submanifold $N$ of $M$, a
  flat section $E$ of $L|_N$, and a subbundle $\mathcal{D}$ of $TM|_N$ such that $\mathcal{D}
  \oplus TN = TM|_N$. Then we can extend $E$ to a neighborhood of $N$ in such
  a way that it satisfies on $N$: $\nabla E =0$ and $\nabla_X \nabla_X E =0$
  for any vector field $X$ of $M$ such that $X|_N$ is a section of $\mathcal{D}$.
  To see this, introduce a coordinate chart $(U, x_i, y_j)$ of $M$ and a unitary frame
  $s:U \rightarrow L$ such that $N \cap U = \{x_1 =
  \ldots = x_k =0 \}$, $ ( \partial_{x_1}, \ldots , \partial_{x_k})$ is a
  frame of $\mathcal{D}$ and $s$ extends $E$. Then the section we are looking for is $e^{i \varphi} s$ with
  $$ \varphi = \sum_{i=1}^{k} \beta_i (0,y) x_i +\tfrac{1}{2} \sum_{i,j =1}^{k}
  (\partial_{x_j} \beta_i )(0,y) x_i x_j + \bigo (|x|^3) $$
  where the $\beta_i$'s are the functions in $\Ci(U)$
  such that $ \nabla_{\partial_{x_i}} s =
  \frac{1}{i} \beta_i \, s$. Applying this to $M^2$, $L \boxtimes \con {L}$ and $\op{diag}M$ instead
  of $M$, $L$, $N$ concludes the proof. 
\end{proof}

\begin{defin} \label{def:heis-semicl-oper}
  A semiclassical Heisenberg pseudodifferential operator of order
  $m \in \R$ is a family of operators
  $(P_k : \Ci ( M, L^k) \rightarrow \Ci ( M , L^k) , \; k \in \N)$ whose
  Schwartz kernels have the form
\begin{gather} \label{eq:def_heis_global}
  F^k (x,y)  \phi(x,y)    K_{k^{-1/2}}(x,y)  + \biginf ( k^{-\infty})
\end{gather}
where $ ( V, F)$ satisfies the conditions of Lemma \ref{lem:existenceVF},
  $\phi \in \Ci_0(V)$ is equal to $1$ on a neighborhood of the diagonal and
  $(K_h , \; h \in (0,1])$ is the Schwartz kernel family of a semiclassical pseudodifferential operator
  $(Q_h) \in \psc^m ( M)$.

  The principal symbol $\si (P)$ of $(P_k)$ is defined as the
  principal symbol of $(Q_h)$. 
\end{defin}

We denote by $\Psi ^{m}_{\op{Heis}} ( L, \nabla )$ the space of semiclassical
Heisenberg pseudodifferential operators of
order $m$. For any $P \in \Psi ^{m}_{\op{Heis}} ( L, \nabla )$, 
for any fixed $k$, $P_k$ is a pseudo-differential operator of order $m$, so
$P_k$ act on $\Ci ( M , L^k)$ and on $\mathcal{C}^{-\infty} ( M, L^k)$.
 The definition clearly does not depend on the choice of the cutoff
function $\phi$. It
neither doesn't depend on the choice of $F$ as will be explained below. To
compare with the twisted pseudodifferential operators, observe first that the
section $F$ in \eqref{eq:def_alt_twisted_pseudo} satisfies a weaker condition
than in Definition  \ref{def:heis-semicl-oper} and second in
\eqref{eq:def_alt_twisted_pseudo}, the
Schwartz kernel of $Q$ is evaluated at $h = k^{-1}$, whereas in 
\eqref{eq:def_heis_global}  we have $h = k^{-1/2}$.

By defining globally the Heisenberg
pseudodifferential operators in terms of scalar pseudodifferential operators as in Definition
\ref{def:heis-semicl-oper} instead of the local oscillatory integrals
\eqref{eq:intro_Def_Heisenberg}, we avoid the usual discussions on the coordinate changes and
the principal symbol and we deduce easily the following three facts: 
\begin{itemize}
  \item[-]   If $P\in \pheis^m (L, \nabla )$, then $\si (P) =0$ if and
    only if $P \in  k^{-\frac{1}{2} } \pheis^{m-1}  
(L, \nabla )$. 
\item[-] $  \bigcap_m  k^{-\frac{m}{2}} \pheis^{-m} (L,\nabla) = k^{-\infty}
  \Psi^{-\infty} (L) $.
\item[-]    If $P\in \pheis^m (L, \nabla)$ and $\rho \in \Ci ( M^2)$ is such that
  $\op{supp} \rho \cap \op{diag} M = \emptyset$, then the kernel $(x,y)
  \rightarrow \rho (x,y)
      P_k (x,y)$ is smooth and negligible.
\end{itemize}
Unfortunately, the definition \ref{def:heis-semicl-oper} does not allow to
deduce the composition properties of the Heisenberg operators from the one of
the semiclassical pseudodifferential operators.

By Lemma \ref{lem:existenceVF}, $F$ is uniquely
defined modulo a factor $e^{i \psi}$ with $   \psi \in \Ci ( U^2)$ vanishing to third order along the
diagonal. Write
\begin{gather} \label{eq:dec_psi}
\psi (x,y)  = \sum_{|\al|=3} \psi_{\al} (x,y) ( x - y)^{\al} 
\end{gather}
with
smooth coefficients $\psi_\al$.
For any symbol $a \in S^{\infty} ( U^2, \R^n)$, let $I(a)$ be the oscillatory
integral
\begin{gather} \label{eq:oscil_int}
I(a) (h, x, y) = \int e^{i h^{-1} \xi \cdot (x-y)} a (
  h,x ,y ,\xi) \; d\xi .
\end{gather}
\begin{lemme} \label{lem:integral_changement_Section}
    For all 
  $a \in S^m ( U^2, \R^n)$, $e^{i h ^{-2} \psi(x,y) } I(a)( h ,x,y) =
  I(b) ( h ,x,y)$ with $b \in S^m (U^2, \R^n)$ having the asymptotic expansion
  \begin{gather} \label{eq:expansion_changement_F}
    b = \sum_{\ell =0 } ^\infty \frac{ h^\ell}{\ell!} L^\ell (a) , \qquad
  \text{ with } \quad L
  = \sum_{|\al|=3} \psi_{\al} (x,y) \partial_{\xi}^\al .
\end{gather}
In particular, if $a \in S^m_{\op{sc}} (U^2 , \R^n)$, then $b \in  S^m_{\op{sc}} (U^2 , \R^n)$.
\end{lemme}

\begin{proof} By integration by part, $(x_i - y_i) I(a) = i hI
  ( \partial_{\xi_i} a)$, so by \eqref{eq:dec_psi}, we have $i h^{-2} \psi I(a) = h I (L( a))$ with $L$ given
  by \eqref{eq:expansion_changement_F}.
  By Taylor formula, we have 
  $$  e^{i h^{-2} \psi} =  \sum_{\ell = 0 } ^{N} \frac{ (i h ^{-2} \psi )^{\ell}
  }{\ell !} + (ih ^{-2} \psi )^{N+1} r_N $$
 with the remainder $$r_N (h,x,y) = \frac{1}{N!} \int_0^1 e^{it h^{-2} \psi (x,y)} (1-t)^{N}
 \; dt.$$
 It comes that
  $$ e^{i h^{-2} \psi} I(a)  = \sum_{\ell = 0 } ^{N} \frac{h ^{\ell}
  }{\ell !} I ( L^\ell (a)) + h ^{N+1} r_N I ( L^{N+1} (a))$$
  We claim that for any $k \in \N$, when $N$ is sufficiently large,
$r_N  I(L^{N+1} (a)) $ is of class
  $\mathcal{C}^k$ and for any $\al \in \N^{2n}$ with $|\al| = k$, $h^{2|\al|} \partial^\al_{x,y} (r_N
  I(L^{N+1} (a)) ) = \bigo (1)$.

  Indeed, $r_N$ is smooth and $h^{2|\al|} \partial^\al_{x,y} r_N =
  \bigo (1)$. Furthermore, since $I(a)$ is a genuine integral for
  $m < -n$, by derivating \eqref{eq:oscil_int} under the integral sign $k$
  times, it comes for   $k + m < -n$ that  $I(a) \in \mathcal{C}^k$ and  $h ^{|\al|} \partial^\al_{x,y}
  I(a) =  O(1)$ for $|\al| = k$. We deduce the claim by using that $L^{N+1}
  (a)  $ is a symbol of order $m - 3 (N+1)$.

 So for any $ b \in S^m (U^2, \R^n)$ having the asymptotic expansion
  \eqref{eq:expansion_changement_F}, we have that $e^{i h^{-2} \psi} I(a) = I(b) + \rho$ with $\rho \in
  h^{\infty} \Ci (U^2)$, and we
  can absorb $\rho$ in $I(b)$ by modifying $b$ by a summand in $h^{\infty}
  S^{-\infty} ( U^2, \R^n)$.  
\end{proof}

{\cc By Lemma \ref{lem:integral_changement_Section}, Definition \ref{def:heis-semicl-oper} is independent of $F$ in the sense that for any $(V, F)$ satisfying the condition
\ref{lem:existenceVF}, for any $P \in  \Psi ^{m}_{\op{Heis}} ( L, \nabla )$,
\eqref{eq:def_heis_global} holds for a convenient $(Q_h)$. Moreover, since $b = a + \bigo (h)$ in \eqref{eq:expansion_changement_F}, the principal symbol map
$$\si : \Psi ^{m}_{\op{Heis}} ( L, \nabla ) \rightarrow 
S_{\op{ph}}^m ( M, T^*M) $$ is
 also independent of  $F$. }

\section{Local expressions} \label{sec:local-expressions}

Let us explain how we recover the local expression
\eqref{eq:intro_Def_Heisenberg} of the introduction. Let
$(U, x_i)$ be a local chart of $M$ and $s \in \Ci ( U, L)$ such
  that $|s|
  =1$ on $U$. Let $\beta \in \Om^1 ( U, \R)$ be the connection one-form, $\nabla s
  = \frac{1}{i} \beta \otimes s$.  Then we easily check that the section $F \in \Ci ( U^2, L \boxtimes
  \con{L})$ given by 
\begin{gather} 
\label{eq:local_F}
 F (x,y) = e^{i \beta \bigl( \tfrac{x+y}{2} \bigr) \cdot (x-y)  } s(x) \otimes \con{s}
 (y)
 \end{gather}
satisfies the condition of Lemma \ref{lem:existenceVF}. Consequently, the
Schwartz kernel of an operator  $P \in \Psi^m_{\op{Heis}} ( L, \nabla)$ has the form $K_k
s^k \boxtimes \con{s}^k$ on $U^2$ with
\begin{gather} \label{eq:loc_racine} 
  K_k (x,y)  = e^{ i k \be \bigl( \tfrac{x+y}{2} \bigr) \cdot (x-y) } \Bigl(
  \frac{\sqrt{k}}{ 2 \pi } \Bigr)^n \int e^{i \sqrt{k} \; \xi \cdot
    (x-y) } a ( k^{-\frac{1}{2}} ,x,y, \xi ) \; d\xi 
\end{gather}
with  $a \in S^m_{\op{sc}} ( U^2, \R^n)$. {\cc The principal symbol of $P$ is
  determined over $U$ by
  $ \si (P) (x, \xi) = a ( h, x, x, \xi) + \bigo (h) .$}

Of course, we can assume that $a$ does
not depend on $y$ (resp. $x$)  or that it is on the Weyl form $a( h, x, y, \xi)
= b ( h, 
\frac{1}{2}(x+y) , \xi)$ with $b \in S^m_{\op{sc}} ( U, \R^n)$. In this last
case, we recover exactly the expression \eqref{eq:intro_Def_Heisenberg}.

{\cc Let us denote by  $\op{Op}_{\op{Heis}} ( a )$  the operator with Schwartz kernel
\eqref{eq:loc_racine} acting on sections of $L^k \rightarrow U$.
Introduce the rescaled covariant derivatives
$$ \tpi_j =\tfrac{1}{i \sqrt k} \nabla_j = \tfrac{1}{i \sqrt k } \partial_j -
\sqrt k \beta_j, \qquad j =1, \ldots, n$$
acting as well on sections of $L^k \rightarrow U$. For any $i,j = 1, \ldots,
n$, let $$\om_{ij}  = \partial_{x_i} \beta_j
- \partial_{x_j} \beta_i \in \Ci (U),$$ so that $\om = d \be = \sum_{i<j}
\om_{ij} dx_i \wedge dx _j$. For any $f \in \Ci(U)$, let $M_f$ denote the
multiplication by $f$, acting on sections of $L^k \rightarrow U$.  
\begin{lemme} \label{lem:corder}
  For any $ a \in  S^m_{\op{sc}} ( U^2, \R^n)$,
  \begin{itemize}
    \item   $ \tilde \pi_j \circ \op{Op}_{\op{Heis}} ( a
  ) = \op{Op}_{\op{Heis}} ( c) $ for some
  $ c \in S^{m+1} _{\op{sc}} ( U^2, \R^n)$
such that
 $$ c( h, x, x,
 \xi) = \Bigl(\xi_j + \tfrac{i}{2} \sum_{\ell =1}^{n} \om_{j \ell}
 \partial_{\xi_{\ell}} \Bigr) a ( h , x, x, \xi )  +  \bigo (h). $$
 \item $M_f \circ  \op{Op}_{\op{Heis}} ( a ) = \op{Op}_{\op{Heis}} ( c)$ with 
   $c(h,x,y) = f (x) a ( h,x,y)$.
 \end{itemize}
\end{lemme}

\begin{proof} Derivating \eqref{eq:loc_racine} with respect to $x_j$, we get first that $\tilde \pi_j \circ \op{Op}_{\op{Heis}} ( a
) = \op{Op}_{\op{Heis}} ( b
)$ with 
$$ b (h,x,y, \xi) = \bigl( h^{-1} \psi_j (x,y) + \xi_j + \tfrac{h}{i}
\partial_{x_j} \bigr) \, a
(h,x,y, \xi)$$
where $$  \psi_j (x,y)  =  \tfrac{1}{2} (\partial_j \be) (\tfrac{1}{2} (x + y) ) \cdot (x - y ) +
  \beta_j (\tfrac{1}{2} (x + y) ) - \beta_j (x).$$
 Taylor expanding along $x = y$, we get
  \begin{gather*}
  \psi_j (x,y) =  \tfrac{1}{2} \sum_\ell \om_{j\ell} (x) ( x_\ell - y_\ell) +
  \sum_{\ell,m} r_{\ell m } (x,y)
( x_\ell - y_\ell)(x_m- y_m) 
\end{gather*}
Integrating by part, it comes that  $\tilde \pi_j \circ \op{Op}_{\op{Heis}} ( a
) = \op{Op}_{\op{Heis}} ( c
) $ with 
\begin{xalignat*}{2} 
c (h,x,y, \xi) &  = (  \xi_j + \tfrac{i}{2}  \sum_\ell \om_{j\ell} ( x)
\partial_{\xi_\ell} \Bigr)  a (h,x,y, \xi)  \\
 & + h  \Bigl(   \tfrac{1}{i} \partial_{x_j}  -  \sum_{\ell, m }
r_{\ell m } (x,y)
\partial_{\xi_\ell} \partial_{\xi_m} \Bigr)  a (h,x,y, \xi) .
\end{xalignat*}
This ends the proof of the first assertion. The second one is obvious.
\end{proof}

Consequences of Lemma \ref{lem:corder} will be drawn in Section
\ref{sec:heis-diff-oper}.
Let us anticipate slightly and derive the equation that the resolvent symbol
has to satisfy. So introduce a Riemannian metric $g$ on $U$ and define the
corresponding Laplacian 
 $$ k^{-1} \Delta_k = \tfrac{1}{2 \sqrt g }  \sum_{j, \ell  =1 }^{n}  \tpi_j
g^{j \ell}
\sqrt g \, \tpi_{\ell} $$
Given $z \in \C$, our goal is to
find a Heisenberg pseudodifferential operator $P$ of order $-2$  such that
\begin{gather} \label{eq:reso}
(k^{-1} \Delta_k  - z) \circ P_k  = \op{id} . 
\end{gather}
By Lemma  \ref{lem:corder}, $ (k^{-1} \Delta_k  - z) \circ P_k $
is a Heisenberg pseudodifferential operator of order $0$. Moreover, $
\op{id} = \op{Op}_{\op{Heis}} ( 1)$ is a Heisenberg pseudodifferential operator of order $0$
as well. The main step in solving \eqref{eq:reso} is to find the principal
symbol $\si$ of $P_k$, that is to solve \eqref{eq:reso} modulo $\Psi^{-1}_{\op{Heis}}
( L, \nabla)$.  
By Lemma   \ref{lem:corder}, the principal symbol of $k^{-1} \Delta_k
P_k$ is $\square \si$ where 
\begin{gather} \label{eq:box}
\square = \tfrac{1}{2} \sum_{j, \ell} g^{j\ell} \Bigl(   \xi_j + \tfrac{i}{2}  \sum_m \om_{jm} 
\partial_{\xi_m} \Bigr) \Bigl(  \xi_\ell + \tfrac{i}{2}  \sum_p \om_{jp}
\partial_{\xi_p} \Bigr) .
\end{gather}
So we are looking for $\si$ such that
\begin{gather} \label{eq:resolvent_symbol_equation}
   ( \square
-z ) \si = 1, \qquad \si \in S^{-2}_{\op{ph}} ( U, \R^n) . 
\end{gather}
If $\om =0$, $\square $ is
merely the multiplication by  $ \tfrac{1}{2} |\xi |^2 = \tfrac{1}{2} \sum g^{j\ell} \xi_j
 \xi_{\ell}$, and we have a solution as soon as $z \notin \R_{\geqslant
   0}$. 
The general case is more complicated because of the derivatives in
\eqref{eq:box}. As we will see, in the case where $\om$ and $g$ are
compatible, \eqref{eq:resolvent_symbol_equation}  has a solution as soon as $z \notin ( \tfrac{d}{2}
 + \N)$. Our strategy to solve \eqref{eq:resolvent_symbol_equation} will be to
 rewrite it as an equality between operators.}

Let us come back to the Heisenberg Schwartz kernel \eqref{eq:loc_racine} and
derive another useful expression. {\cc Assume that  $a$ is on the Weyl form, that is $a( h, x, y, \xi)
= b ( h,
\frac{1}{2}(x+y) , \xi)$, then we have that
\begin{gather} \label{eq:Weyl_form}
K_k (x,y) =   \Bigl(
  \frac{k}{ 2 \pi } \Bigr)^n \int e^{i k \; \xi \cdot
    (x-y) }  \tilde{b} \bigl( k^{-1} , \tfrac{1}{2} ( x+ y) , \xi \bigr) \; d \xi
\end{gather}
  where $\tilde {b} ( h, x, \xi ) = b ( \sqrt h, x, h^{-\frac{1}{2}} ( \xi -
  \be (x)))$
\begin{proof}[Proof of Formula \eqref{eq:Weyl_form}] By the change of variable
  $\xi \rightarrow \sqrt k \, \xi$ in \eqref{eq:loc_racine}, we get
$$  K_k (x,y)    = e^{ i k \be \bigl( \tfrac{x+y}{2} \bigr) \cdot (x-y) } \Bigl(
  \frac{k}{ 2 \pi } \Bigr)^n \int e^{i k \; \xi \cdot
    (x-y) } a ( k^{-\frac{1}{2}} ,x,y, \sqrt{k} \; \xi ) \; d\xi $$
Then, by the change of variable $\xi \rightarrow \xi - \beta ( \frac{x+y}{2})
$, we obtain
\begin{xalignat*}{2} 
 K_k (x,y)   &  =   \Bigl(
  \frac{k}{ 2 \pi } \Bigr)^n \int e^{i k \; \xi \cdot
    (x-y) } a \bigl( k^{-\frac{1}{2}} ,x,y, \sqrt{k} \bigl(   \xi  - \be \bigl(
  \tfrac{x+y}{2} \bigr) \bigr)  \bigr)\; d\xi \\
 & =  \Bigl(
  \frac{k}{ 2 \pi } \Bigr)^n \int e^{i k \; \xi \cdot
    (x-y) }  b ( k^{-\frac{1}{2}} , \tfrac{x+y}{2}, \sqrt k ( \xi -
  \be (\tfrac{x+y}{2}))) \; d \xi
\end{xalignat*}
and we recognise  \eqref{eq:Weyl_form}.
\end{proof}


The right-hand side of \eqref{eq:Weyl_form} is the Schwartz kernel of   a semiclassical
pseudodifferential operator at $k = h^{-1}$ with a Weyl symbol
$\tilde{b}$. For this reason we  call $\tilde{b}$ the {\em effective} symbol.
Unfortunately $\tilde{b}$ does not satisfy the usual symbolic estimates but
some exotic ones.}

 Let us introduce the symbol semi-norms of $S^m (U, \R^n)$, 
$$ \| a \|_{m, \ell, K } = \max _{ |\al | + | \be | \leqslant \ell}
\sup_{x\in K , \xi \in \R^n} |  \partial_x^{\al} \partial_{\xi}^{\be} a
(x,\xi) | \langle \xi \rangle^{-m  + | \be |}
$$
where $K$ is a compact subset of $U$. 
\begin{lemme} \label{lem:exotic}
  For any $m \in \R$, $\al, \be \in \N^n$ and compact subset
  $K$ of $U$, there exists $C>0$ such that for any $a   \in \Ci (U \times \R^n)$, the function $\tilde{a}
  (h, x, \xi ) = a(  x, h^{-\frac{1}{2}} ( \xi - \be (x) )$ satisfies 
$$ | \partial_x^{\al} \partial_{\xi}^{\be} \tilde{a} (h, x, \xi) | \leqslant
C
\| a \|_{m, \ell,K} h^{-\frac{1}{2} ( m_+ + \ell )}  \langle \xi
\rangle^{m - | \be |}  $$
for all $ 0 < h \leqslant 1$, $x \in K$, $\xi \in \R^n $
with $\ell = | \al | + | \be|$, $m_+ = \max ( m,0)$. 
\end{lemme}
\begin{proof} For any $0 < \ep \leqslant 1 $, we have $\langle
  \eta \rangle \leqslant \langle \ep^{-1}  \eta \rangle \leqslant \ep^{-1}
  \langle \eta \rangle$. Furthermore, if $x \in K$, $ C^{-1} \langle \xi
  \rangle \leqslant \langle \xi - \be (x) \rangle \leqslant C \langle \xi
  \rangle$. So for any $m \in \R$, 
\begin{gather} \label{eq:e1}
  \langle \ep^{-1} ( \xi - \be (x) ) \rangle^m \leqslant C_m \ep^{- m_+}
\langle \xi \rangle ^m . 
\end{gather}
The derivatives of $ \tilde{a} _{\ep}( x, \xi ) = a ( x, \ep^{-1} (\xi - \be
(x)))$  have the form  $ \partial_x^\al \partial_\xi^\be \tilde a_{\ep}  = \tilde b_{\ep} $
with 
\begin{gather} \label{eq:e2}
b = \sum _{\al', \be'} \ep^{- | \be'|} f_{\al', \be' } \partial_x^{\al'}
\partial_\xi^{\be'} a  
\end{gather}
where the coefficients $f_{\al',\be'}$ are in $\Ci (U)$ and don't depend on $a$, and we sum over
the multi-indices satisfying $ \be \leqslant \be'$ and $|\al' | + | \be'|
\leqslant | \al | + | \be|$.  So for $x \in K$, 
\begin{xalignat*}{2} 
|  \partial_x^\al \partial_\xi^\be \tilde a_{\ep} (x, \xi) | & \leqslant
C \sum_{\al', \be' } \ep^{- | \be'|} \| a \|_{m, \ell, K} \langle \ep^{-1} (
\xi - \be (x)) \rangle^{m - |\be'|}  \quad \text{ by }  \eqref{eq:e2} \\
&  \leqslant
C \sum_{\al', \be' } \ep^{- | \be'|} \| a \|_{m, \ell, K} \langle 
\xi \rangle^{m - |\be'|} \ep^{-m_+}  \quad \text{ by }  \eqref{eq:e1} \\
&  \leqslant
C \ep^{ - ( | \al | + | \be| )}  \| a \|_{m, \ell, K}  \langle  \xi \rangle^{m - |\be |} \ep^{-m_+}
\end{xalignat*}
 because $|\be | \leqslant | \be'| \leqslant |\al | + | \be|$
and we conclude by setting $\epsilon = h^{\frac{1}{2}}$.
\end{proof}

So when $b $ belongs to $S^{m} (U, \R^n)$, the semi-norms $\| b \|_{m, \ell,
  K}$ are finite and by Lemma \ref{lem:exotic}, $\tilde{b}$ belongs to the class $S_{\delta}$ with exponent $\delta = \frac{1}{2}$, that is at each derivative we loose a factor $h^{-\delta}$. 
Recall that $\delta =\frac{1}{2}$ is the critical exponent: the space of semiclassical pseudodifferential operators with symbol in
$S_{\delta}$ is an algebra for $\delta \in [0, \frac{1}{2} ]$, but the standard
asymptotic expansions of the symbolic calculus only hold for $ \delta \in
[0,\frac{1}{2} ) $, cf. for instance \cite[Proposition 7.7]{DiSj}. As we will see in
Section \ref{sec:heis-diff-oper} and in \cite{oim_Heis}, the Heisenberg
pseudodifferential operators form an algebra and have an associated symbol
calculus, but this can not be deduced from \eqref{eq:Weyl_form} and the usual composition rules of
pseudodifferential operators. Nevertheless, Formula \eqref{eq:Weyl_form} has some useful
consequences, as we will see in the next section.

\section{Mapping properties} \label{sec:mapping-properties}

 Recall the
definition \eqref{eq:l2norm} of the $L^2$-norm with a volume element
independent of $k$.

\begin{theo} \label{theo:mapping_property} 
  For any $P \in \Psi^0_{\op{Heis}} ( L, \nabla)$, there exists $C>0$ such that for
  any $k$,  $ \| P_k \| _{\mathcal{L} ( L^2
    ( M , L^k))} \leqslant C $.
\end{theo}

\begin{proof}
Introduce a finite atlas $(U_i, \phi_i)$ of $M$ with functions $\varphi_i,
\psi_i \in \Ci_0( U_i) $ such that $\sum \varphi_i =1$ and $\op{supp} \varphi_i \subset \op{int}
\{ \psi_i =1 \}$. Write 
\begin{gather} \label{eq:decomp_P_local}
P = \sum \psi_i P \varphi_i + Q.
\end{gather}
Since $\sum \psi_i
(x) \varphi_i (y) =1$ when $x$ is close to $y$, $Q$ is in $k^{-\infty}
\Psi^{-\infty}$. Identifying $U_i$ with $\phi_i (U_i)$,  the Schwartz kernel
of $\psi_i P \varphi_i $ has the form \eqref{eq:Weyl_form} with a symbol $\tilde{b}_i $
satisfying by Lemma \ref{lem:exotic}
$$  | \partial^\al_x \partial^\be_\xi \tilde b_i (h, x, \xi ) |
\leqslant h^{-\frac{1}{2} ( |\al| + |\be|)} C_{\al, \be}.$$
{\cc By Calderon-Vaillancourt for semiclassical pseudodifferential operators, cf.
as instance} \cite[Theorem 7.11]{DiSj},  $\psi_i P \varphi_i= \bigo (1) : L^2 ( \R^n)
\rightarrow L^2 ( \R^n)$. 
\end{proof}

Another consequence of Expression \eqref{eq:Weyl_form} and Lemma
\ref{lem:exotic} is the following important fact that we will need later. 

\begin{lemme} \label{lem:residuel}
  $k^{- \infty} \Psi^{-\infty}(L) $ is a bilateral ideal of
  $\pheis ( L, \nabla)$.
\end{lemme}
\begin{proof}  
  Consider a pseudodifferential operator $A(h)$ of $\R^n$ with the Schwartz  kernel
  $(2 \pi h )^{-n} \int e^{ih^{-1} \xi (x-y) }
a(h,x,y,\xi) d\xi $ where the amplitude $a(h,x,y, \xi) $ is zero if $|x|+|y| \geqslant
C $ and satisfies
\begin{gather} \label{eq:-partialal_x-y}
  | \partial^\al_{x,y} a(h,x,y,\xi) |  \leqslant h^{-|\al|} C_{\al}
\langle \xi \rangle ^m , \qquad \forall \, \al. 
\end{gather}
Then, with the usual regularisation of oscillatory integrals by
integration by part, one
proves that {\cc for any $\al \in \N^n$, there exists $C_\al'$ such that for any $u
\in \Ci_0 ( \R^n)$ and $h \in (0,1]$,   
$$ h^{|\al|} \sup_{x \in \R^n}  | \partial^{\al}_x (A(h) u (x)) |   \leqslant
C'_{\al} \max _{| \be| \leqslant m + n +1 + | \al|} \sup_{x \in \R^n}  h^{|\be|}
|\partial_x^\be u(x) | .$$ }
So for any family of operators $B(h) : \Ci ( \R^n) \rightarrow \Ci ( \R^n)$,
$h \in (0,1]$, 
if $B(h)$ has a compactly supported smooth kernel in $\biginf ( h^{\infty})$, then the
same holds for $A (h) \circ B(h)$.

{\cc This implies that $k^{- \infty} \Psi^{-\infty}(L) $ is a left ideal of
  $\pheis ( L, \nabla)$. Indeed, any Heisenberg pseudodifferential operator
  $P$ may be decomposed as in \eqref{eq:decomp_P_local}, where each  $\psi_i P
  \varphi_i$ is a semiclassical pseudodifferential operator with
  a symbol satisfying \eqref{eq:-partialal_x-y} by Lemma \ref{lem:exotic}. To
  prove that $k^{- \infty} \Psi^{-\infty}(L) $ is a right ideal, merely take formal adjoints. }
\end{proof}



To end this section, let us extend the mapping property to the Sobolev space.
We denote by $\| \cdot \|_m$ the $m$-th semiclassical Sobolev norm of sections of $L^k$,
defined as in Section \ref{sec:twist-pseud-oper}. 

\begin{theo} \label{theo:map_Sobolev}
  For any $m , p \in \R$ and any ${\cc P} \in \pheis^m ( L, \nabla)$, there exists $C>0$
  such that for any $k$, 
  $$ \| P_k u  \|_{p}   \leqslant C k^{\frac{1}{2} m_+} \| u
  \|_{p+m}  , \qquad \forall u \in \Ci ( M , L^k).$$
\end{theo}

Since for any $k$, $P_k$ is  a pseudodifferential operator of order
 $m$ of $L^k$, we already know that $P_k$ is continuous $H^{p+m} ( M , L^k) \rightarrow H^{m} (M, L^k)$.
 Theorem \ref{theo:map_Sobolev} gives a uniform estimate with respect to $k$.
 
\begin{proof} It suffices to prove that for any $E \in \ptw ^{p - m } (L)$
  and  $E' \in \ptw^{-p} ( L)$ one has
  \begin{gather} \label{eq:toprove0}
    E' _k P_k E_k = \bigo (
 k^{\frac{1}{2} m_+}) : L^2(M, L^k)
 \rightarrow L^2 (M , L^k) .
\end{gather}
For this it suffices to prove that for any chart domain $U$ of $M$
 and functions $\rho _j \in \Ci _0 ( U)$, $j=1,...,4$, one has
\begin{gather} \label{eq:toprove}
 \rho_1  E' _k \, \rho_2 P_k \, \rho_3 E_k \, \rho_4  = \bigo (
 k^{\frac{1}{2} m_+}) : L^2(M, L^k)
 \rightarrow L^2(M, L^k) .
\end{gather}
To show that \eqref{eq:toprove} implies \eqref{eq:toprove0}, write $P$ on the form \eqref{eq:decomp_P_local},  $E' \psi_i =
 \tilde{\psi}_i E' \psi_i +  ( 1 - \tilde{\psi}_i ) E' \psi_i$ with $\tilde{\psi}_i \in \Ci_0 ( U_i)$ such that $\op{supp} \psi_i \subset \op{int} \{
 \tilde{\psi}_i =1 \}$ and similarly for $E$, and use that $k^{-\infty}
 \Psi^{-\infty} (L)$ is an ideal of both  $\pheis^{\infty}(L, \nabla) $ and
 $\ptw^{\infty} (L)$. 

 As in \cite[Definition 7.5]{DiSj}, for $\delta \in [0,1]$ and $m : \R^n
 \rightarrow [0, \infty )$ an order function, let $S_\delta (m)$ be the space
 of families $(a(h), \;  h \in (0,1])$ of $\Ci (\R^n)$ such that $ | \partial^\al a(h,x) |
 \leqslant C_{\al} h^{-\delta |\al|}  m(x)$.  {\cc Recall that $m$ is an order
 function means that for some positive constants $C$, $N$, we have $m (x)
 \leqslant C \langle x -y \rangle^N m (y)$ for any $x, y\in \R^n$. The
 order functions we need are the functions $m_r$ defined by 
$$ m_r : \R^{2n} \rightarrow [0, \infty ), \qquad m_r(x,\xi) = \langle \xi
 \rangle^r $$ 
for $r \in \R$.}

Identify $U$ with $\phi (U)$ and denote by $\op{Op}_k ( \tilde{b} ) $ the
operator with kernel \eqref{eq:Weyl_form}.
By Lemma \ref{lem:exotic}, for any $\rho, \rho ' \in
\Ci_0 ( U)$,  $ \rho  E' _k \,
\rho'$,  $k^{-\frac{1}{2} m_+} \rho  P_k \, \rho'$ and $\rho E_k \, \rho'$ are
equal to $\op{Op}_k ( \tilde{b})$  with $\tilde{b} $ in $S_{0} ( \langle \xi \rangle^{-p}
)$, $ S_{1/2} ( \langle \xi \rangle^m)$ and $S_{0} ( \langle \xi
\rangle^{p-m})$ respectively. By \cite[Proposition 7.7, Theorem 7.9]{DiSj},
their product is equal to $\op{Op}_k (c)$ with  $c \in  S_{1/2} ( 1)$, which
proves \eqref{eq:toprove} by \cite[Theorem 7.11]{DiSj}.
\end{proof}

Actually, Theorem \ref{theo:map_Sobolev} can be improved if we use Sobolev norms associated
to the covariant derivative $\nabla$ instead of the semiclassical Sobolev norms. For
 instance, for any $m \in \N$, any $Q \in \pheis^{-m} ( L,\nabla)$ and any vector
 fields $X_1$, \ldots, $X_m$ of $M$, we will see in Proposition
 \ref{prop:heis-diff-pseudo}, that $P= (k^{-m/2} \nabla_{X_1} \ldots  \nabla_{X_m}
 Q_k)$ belongs to $\pheis^0 ( L, \nabla)$, so by Theorem
 \ref{theo:mapping_property},
\begin{gather} \label{eq:mieux}
 k^{-m/2} \nabla_{X_1} \ldots  \nabla_{X_m}
 Q_k = \bigo (1) :L^2 ( M , L^k) \rightarrow L^2 ( M , L^k) 
\end{gather}
To compare, Theorem \ref{theo:map_Sobolev} only implies that the norm of $P_k$
in $\mathcal{L} (L^2 ( M , L^k))$ is in $\bigo ( k^{m/2})$. The generalisation
of \eqref{eq:mieux} to fractional exponents $m$ not necessarily nonnegative
will be given in \cite{oim_Heis}.

\section{A product associated to an antisymmetric bilinear form} \label{sec:isotropic-algebra}

Let $E$ be a $n$-dimensional real vector space and $A \in \wedge^2
E^*$. Later, we will choose $E =T_xM$ with $A = \om (x)$. Introduce the covariant
derivative of $E$ 
\begin{gather} 
\nabla^A = d + \tfrac{1}{i} \beta, \qquad \text{ where } \be \in \Om^1 ( E,
\R), \quad  \beta (x) (Y) = \tfrac{1}{2} A (x, Y).
\end{gather}
{\cc So for any $Y \in E$ considered as a constant vector field of $E$,
the covariant derivative $\nabla_Y^A$ acts on $\Ci (E)$ by $Y + \frac{1}{i}
\be ( \cdot ) (Y)$.}   
The curvature of $\nabla ^A$ is
$\frac{1}{i} A$, that is
\begin{gather} [\nabla_X^A, \nabla_Y^A ] = \tfrac{1}{i} A(X,Y), \qquad \forall
  \, X, Y
  \in E 
\end{gather}
{\cc Indeed, 
  $[X,Y ] = 0$ and
  the de Rham derivative of $\be$ is $A$, viewed as a constant $2$-form.}

We will define for any tempered distribution $g \in \mathcal{S}'(E^*)$ an
operator $g ( \frac{1}{i} \nabla^A)$.
We assume first that $E = \R^n$. 
For any $g \in \mathcal{S}'( \R^n)$, we denote by $\wh
 g $ and $g^{ \vee}$ its Fourier transform and inverse Fourier transform,
 with the normalisation
 $$ \wh{g} (\xi) =
 \int_{\R^n}
e^{-ix \cdot\xi} g ( x) \; dx, \qquad g ^{\vee} ( x ) = (2 \pi)^{-n} \, \wh{g} (-
x). $$
Let $g ( \tfrac{1}{i} \partial) $ be the operator from
$\mathcal{S} ( \R^n )$ to $ \mathcal{S}'(\R^n)$ such that $g ( \tfrac{1}{i} \partial )  u = v $ if and
only if $  g ( \xi ) \wh u ( \xi ) =  \wh v ( \xi)$.  The Schwartz kernel of $ g ( \tfrac{1}{i} \partial)$ is $g^{\vee}
( x-y)$.

Then for any antisymmetric bilinear form $A$ of $\R^n$, define $g (
\frac{1}{i} \nabla^A )$ as the operator with Schwartz kernel
\begin{gather} 
K_g (x,y) = e^{
  -\frac{i}{2} A ( x,y) } g^{\vee} (x-y).
\end{gather}
Since $g^{\vee} (x-y)$ is a
tempered distribution of $\R^n_x \times \R^n_y$, the same holds for $K_g $, so $ g ( \tfrac{1}{i} \nabla^A)$ is
continuous from $\mathcal{S} ( \R^n )$ to $ \mathcal{S}'(\R^n)$.
We claim that this definition has an intrinsic meaning for $A \in \wedge^2
E^*$ if we consider that $ g
\in \mathcal{S}' (E^*)$ and $g (\frac{1}{i} \nabla^A)$ is an operator $\mathcal{S} ( E) \rightarrow
\mathcal{S}' (E)$. One way to see this is to write for $g$ and $u$ in
$\mathcal{S} ( \R^n)$ 
\begin{gather} \label{eq:def_gnabla}
(g ( \tfrac{1}{i} \nabla^A) u) (x) = ( 2 \pi )^{-n} \int_{\R^n \times \R^n}
e^{-\frac{i}{2} A(x,y) +  i  \xi \cdot (x-y) } g (\xi) u ( y) \;  d y \, d\xi 
\end{gather} 
and to notice that the product $ \xi \cdot ( x-y)$ is well-defined for $\xi
\in E^*$, $x, y \in E$, and 
the measure $dy \, d \xi$ can be interpreted as the canonical volume form of $
E \times E^*  \simeq \R^n_y \times \R^n_\xi$.

Assume again that $E \simeq \R_x^n$, $E^* \simeq \R_{\xi} ^n$  and let
\begin{gather} \nabla_j^A := \nabla_{\partial_{x_j}}^A =
\partial_{x_j} + \tfrac{1}{2i}  \textstyle{\sum_k}
 x_k A_{kj}
\end{gather}
where $(A_{ij} )$ is the matrix of $A$, so $A_{ij} = A(e_i, e_j)$ with
$(e_i)$ the canonical basis of $\R^n$.  

\begin{lemme}  \label{lem:facile_et_essentiel}
For any   $ f \in \mathcal{S}'(E^*)$, we have 
$\frac{1}{i} \nabla_{j}^A \circ f ( \frac{1}{i} \nabla^A ) =  (\xi_j \sharp_A f)   ( \frac{1}{i}
\nabla^A )$ where 
\begin{gather} \label{eq:EEEquation}
\xi_j \sharp_A f   = (\xi_j + \tfrac{i}{2 } \textstyle{\sum_k} A_{jk} \partial_{\xi_k} )
f .
\end{gather}
\end{lemme}


\begin{proof} 
Simply use the identity 
$$ \tfrac{1}{i} (\partial_{x_j} + \tfrac{1}{2i}  \textstyle{\sum_k}
 x_k A_{kj} ) e^{- \frac{i}{2} A (x,y) + i  \xi \cdot ( x-y ) } = ( \xi_j +
 \frac{1}{2i} \sum_k A_{jk} \partial_{\xi_k} )  e^{- \frac{i}{2} A (x,y) +i  \xi
  \cdot ( x-y ) } $$
in  \eqref{eq:def_gnabla} and integrate by part with respect to the variables $\xi_k$.  
\end{proof}

The reason for the notation $g ( \tfrac{1}{i} \nabla^A)$ is that when $g$ is a
monomial, $g ( \tfrac{1}{i} \nabla^A)$ is merely a symmetrization
of covariant derivatives. The precise result is the  following proposition
which is not really needed in the sequel. Notice first that for $g \equiv 1$,   $
g ( \frac{1}{i} \nabla^A ) = \op{id}$ as a direct consequence of the
definition. 
\begin{prop} \label{prop:symmet}
  For any $N \geqslant 1 $ and $X_1, \ldots X_N \in E$, if  $ g = \prod_{i =1 }^{N} f_i$ with $f_i (\xi) = \xi
  \cdot X_i$,  then
  $$g (  \tfrac{1}{i} \nabla^A ) = \frac{(-i)^N}{N!} \sum_{\si \in
    \mathfrak{S}_N} \nabla_{X_{\si(1)}}^A \ldots  \nabla_{X_{\si(N)}}^A ,$$
  where $\mathfrak{S}_N$ is the group of permutations of $1, \ldots, N$. 
\end{prop}

\begin{proof} For any $ X \in E$ we have by
\eqref{eq:EEEquation} with  $f (\xi ) =
\xi \cdot X$ that
\begin{gather} \label{eq:intEEE}
\tfrac{1}{i} \nabla^A_X \circ g (
\tfrac{1}{i} \nabla^A) = ( f \sharp_A \, g )( \tfrac{1}{i} \nabla^A)
\end{gather}
   where $f \sharp_A \, g = ( f  + \tfrac{i}{2} A( X, \partial_\xi ) ) g$. 
Choosing $ g =1$, we obtain the result for $N=1$. We now proceed by induction
over $N$ and 
assume the
 result holds for $N-1$ with $N \geqslant 2$. Thus
 $$ \frac{(-i)^N}{N!} \sum_{\si \in
    \mathfrak{S}_N} \nabla_{X_{\si(1)}}^A \ldots  \nabla_{X_{\si(N)}}^A  = \frac{1}{N}
  \sum_{j=1}^{N}  \tfrac{1}{i} \nabla^A_{X_j} \circ g_j ( \tfrac{1}{i} \nabla^A_X )$$
  with $g_j = g / f_j$. By \eqref{eq:intEEE}, $f \sharp_A
  g_j = f g_j + \frac{i}{2} \sum_{\ell \neq j} A( X, X_\ell) g_{j \ell}$ where
  $g_{j \ell} = g / (f_j f_{\ell})$.  So we have
   $$ \frac{(-i)^N}{N!} \sum_{\si \in
    S_N} \nabla_{X_{\si(1)}}^A \ldots  \nabla_{X_{\si(N)}}^A  = g (
  \tfrac{1}{i} \nabla^A) + \frac{1}{N}
  \sum_{j \neq \ell } A( X_j, X_\ell)  g_{j \ell} ( \tfrac{1}{i} \nabla^A )$$
  and the sum in the right-hand side is zero because $A$ is antisymmetric whereas
  $g_{j \ell} = g _{\ell j}$.
\end{proof}

Let $\dis ^{\infty} (A)$ be the filtered algebra generated by the covariant derivatives
$\nabla_X^A$ where $X \in E$. More explicitly, 
$\dis^{\infty} (A)  = \cup_{m \in \N} \dis^m (A)$ with 
$$\dis^m (A)  = \op{Span} ( \nabla^A_{X_1} \ldots
\nabla^A_{X_\ell} / \; 0 \leqslant \ell \leqslant m,  X_1, \ldots X_{\ell} \in E
) .$$
Let $\C_{\leqslant m} [E^*]$ be the space of complex polynomial functions
of $E^*$ with degree less than or equal to $m$. By Lemma
\ref{lem:facile_et_essentiel} and Proposition \ref{prop:symmet},
\begin{gather} \label{eq:dism-a-=}
\dis^m (A) = \bigl\{ f ( \tfrac{1}{i} \nabla^A ) , \; f
   \in \C_{\leqslant m} [E^*] \bigr\} .
 \end{gather}
By Lemma  \ref{lem:facile_et_essentiel} again, the left
composition by any element of $\dis^{\infty} (A) $ preserves $ \{ g ( \frac{1}{i} \nabla),
\; g \in \mathcal{S}'(E^*) \}$. This defines the product
\begin{gather} \label{eq:prod}
 \sharp_A:  \C [E^*] \times \mathcal{S}'(E^*) \rightarrow \mathcal{S}'(E^*),
\qquad (f \sharp_A g)  ( \tfrac{1}{i} \nabla )  = f ( \tfrac{1}{i} \nabla ) \circ
g ( \tfrac{1}{i} \nabla) 
\end{gather}

In the sequel we will use the basis $( \nabla^{\al} , | \al | \leqslant m )
$ of $\dis^m (A)$, defined by $\nabla^{\al}:= (\nabla_1^A)^{\al
  (1)} \ldots (\nabla_n^A) ^{\al(n)}$, $\al \in \N^{n}$. Clearly
\begin{gather} \label{eq:notationpuissancedieze}
i^{-| \al|} \nabla^\al = f ( \tfrac{1}{i}
\nabla) \qquad \text{ with }f = \xi^{\sharp \al} :=\xi_1^{\sharp \al (1)} \sharp \ldots
\sharp \,\xi_n^{\sharp \al (n)}  ,
\end{gather}
where we have not written the $A$ dependence to lighten the notations. 
Furthermore, if $| \ga| = m$, 
\begin{gather} \label{eq:comp_base}
\xi^{\sharp \ga}  \sharp_A f = \xi^{\gamma}f +  \sum_{|\al| + |\be|  \leqslant
  m, \; | \al | \leqslant m-1} a_{\al,
  \be, \ga} \xi^{\al} \partial_{\xi}^\be f 
\end{gather}
where the coefficients $a_{\al, \be, \ga} \in \C$ depends smoothly (even
polynomially) on $A$, which follows from Lemma \ref{lem:facile_et_essentiel}
again. Actually there is a closed formula for $\sharp_A$, cf.
\eqref{eq:exp_product_A}, but \eqref{eq:comp_base} is enough for our purpose. 

Introduce the space $\pis^m (A) := \{ f ( \frac{1}{i} \nabla) , \; f \in S^m_{\op{ph}}
( E^*) \}$. We have
$$ \dis^m(A) \subset \pis^m (A) , \qquad \dis^m(A)  \circ \pis^p
(A) \subset \pis^{m+ p}(A) ,$$
the second assertion being a consequence of \eqref{eq:comp_base}.
This is all what we need to define in the next section the symbolic calculus
corresponding  to the
composition of differential Heisenberg operators with Heisenberg
pseudodifferential operators. In the case where $A =0$, $\sharp_A $ is the usual
pointwise product of functions. 
In  Lemma \ref{lem:Weyl_algebra}, we will see that when $A$ is nondegenerate so that $n = 2 d$, $\pis^{\infty} (A)$ is an algebra isomorphic to
the Weyl algebra of $\R^{2d}$.

In the companion paper
\cite{oim_Heis}, we will prove that for any $A$,  $\pis^{\infty} (A)$ is a filtered algebra, that
is $\pis ^m(A)  \circ \pis^p(A) \subset \pis^{m+p} (A)$. Moreover
\begin{gather} \label{eq:exp_product_A}
(f \sharp_A \, g )(\xi) = \Bigl[ e^{\frac{i}{2} A ( \partial_{\xi}, \partial_{\eta} ) }
f ( \xi ) g ( \eta) \Bigl]_{\xi = \eta} 
\end{gather}
So $\pis^{\infty} (A)$ is isomorphic with the algebra called the $A$-isotropic
  algebra in \cite[Chapter 4, section 2]{EM}.

Recall the standard and Weyl quantization maps which associate to any $f \in \mathcal{S}' ( \R^{2n})$ the operators $ f( x, \frac{1}{i}
\partial)$ and $f^w ( x, \frac{1}{i} \partial)$  with Schwartz
kernels
$$  ( 2 \pi)^{-n} \int e^{i \xi \cdot ( x-y)} f(x, \xi) \; d \xi  \quad \text{
  and } \quad( 2
\pi)^{-n} \int e^{i \xi \cdot ( x-y)} f(\tfrac{1}{2} (x+ y) , \xi) \; d \xi$$
respectively. {\cc In general, $f ( x, \frac{1}{i} \partial ) $ and $f^{w} (x
  \frac{1}{i} \partial )$ are different, but for the operators we are
  interested in, they coincide.}

\begin{lemme} \label{lem:Weyl_symbol}
  For any $g \in \mathcal{S}' (\R^n)$, we have
  $$ g ( \tfrac{1}{i} \nabla^A ) = f ( x,
\tfrac{1}{i} \partial) = f^{w} ( x, \tfrac{1}{i} \partial )
$$  where $ f ( x, \xi ) = g ( \xi -
\beta(x))$ and $\be (x)$ is defined in \eqref{eq:def_gnabla}, equivalently  $f ( x, \xi ) = g( \xi_1  - \frac{1}{2} A ( x , e_1),
\ldots , \xi_n - \frac{1}{2} A ( x, e_n))$. 
\end{lemme}

\begin{proof}
By the change of variable $\xi \rightarrow \xi + \be (x)$, 
$$ \int e^{i \xi \cdot ( x-y) } g ( \xi - \beta (x)) \; d\xi   =
 e^{ i \be (x) (x-y) } \int e^{i \xi \cdot ( x-y) } g ( \xi ) \;
 d \xi .$$
$A$ being antisymmetric, $\be ( x) ( x-y) =  - \frac{1}{2} A( x,y)$, which
proves that $g ( \tfrac{1}{i} \nabla^A ) = f ( x,
\tfrac{1}{i} \partial)$. The same proof by using this time that  $ \be (\frac{1}{2} (x+ y) ) (x-y) =
- \frac{1}{2} A ( x, y) $ shows that $ g ( \tfrac{1}{i} \nabla^A ) = f^w ( x,
\tfrac{1}{i} \partial)$.
\end{proof}

\section{Heisenberg differential operators} \label{sec:heis-diff-oper}

The algebra $\dheis^{\infty} (L, \nabla ) $ of Heisenberg differential operators consists of families
 of differential
 operators
\begin{gather} \label{eq:family}
 P = (P_k : \Ci ( M , L^k) \rightarrow \Ci ( M , L^k), \; k \in \N
),
\end{gather}
satisfying some conditions given below. It includes the  multiplications by any $f$ in $\Ci (M)$, the normalised
covariant derivatives $k^{-1/2} \nabla_X$ where $X$ is any vector field of $M$ and the multiplication by
$k ^{-1/2}$. It is actually generated by these operators but it will be easier to use the following definition. 

For any $m \in \N$, $\dheis^m (L, \nabla )$ consists of the  families $P$  of differential operators
of the form \eqref{eq:family} such that for
any coordinate chart $(U, x_i)$ and frame $s \in \Ci ( U, L)$ with $|s|=1$, we
have on $U$,
\begin{gather} \label{eq:def_opdif_heis} 
 P_k = \sum_{\substack{\ell \in \N, \; \al \in \N^n \\ \ell + |\al| \leqslant m}} k^{-
  \frac{\ell}{2}} f_{\ell, \al } \tpi^\al
\end{gather}
where $f_{\ell, \al} \in \Ci (U)$, $\tpi^{\al} = \tpi_1^{\al (1)} \ldots \tpi_n^{\al (n)}$ and 
\begin{gather} \label{eq:moment_rescale}
\tpi_i=\tfrac{1}{i \sqrt k} \nabla_i = \tfrac{1}{i \sqrt k } \partial_i -
\sqrt k \beta_i \qquad \text{ with } \quad \nabla s = \tfrac{1}{i} \sum \beta_i
dx_i  \otimes s
\end{gather}
Set
$$\dheis^{\infty} (L , \nabla)  = \bigcup_{m \in \N} \dheis^m (L , \nabla).$$ 
In the sequel to lighten the notations, we omit $(L, \nabla)$ and write
$\dheis^m$, $\dheis^{\infty}$. 
Since $[\tpi_i, \tpi_j ] = \frac{1}{i} ( \partial_i \be_j - \partial_j
\beta_i)$ and $[ \tpi_i, f ] = \frac{1}{i \sqrt k } \partial_i f$, we see that
$$ \dheis^m   \circ \dheis^p   \subset \dheis^{m+p}  .$$
Notice that $\dheis^{\infty} $ has two filtrations: one ascending
$\dheis^m  \subset
\dheis^{m+1} $ and the other descending $k^{-\ell/2} \dheis $, $\ell \in \N$.
The generators $f$, $k^{-1/2} \nabla_X$ and $k^{-1/2}$ have orders $0$,
$1$, $ 1$ for the former and $0$, $0$, $1$ for the latter. 

By the next proposition, $\dheis^{\infty}  $ is contained in
$\pheis^{\infty} ( L , \nabla)  $ and acts on it.
Being Heisenberg pseudodifferential operators, the elements of $\dheis ^{\infty}$  have
a principal symbol, cf Definition \ref{def:heis-semicl-oper}. As we will see, the product of symbols is the fiberwise product
$\sharp$ of
$T^*M$ defined from $\om$. Precisely, we denote by $\sharp_x$ 
the product
\begin{gather} \label{eq:sharp_produit}
\C _{\leqslant m } [T^*_xM] \times S^p
(T^*_xM) \rightarrow S^{m+p} (T^*_xM)
\end{gather}
associated to $\om(x) \in \wedge^2
T^*_xM$ defined in \eqref{eq:prod} and {\cd{with the notation $\C_{\leqslant m } $
introduced before \eqref{eq:dism-a-=}.}}
We will need as well the polynomials $\xi^{\sharp_x \al}$
defined in \eqref{eq:notationpuissancedieze}. 
\begin{prop}  \label{prop:heis-diff-pseudo} $ $ 
 
\begin{itemize} 
\item[-] for any $m \in \N$, $\dheis^m \subset \pheis^m$
\item[-] the principal symbols of the operators of $\dheis^m$ are the functions $f \in \Ci (T^*M)$ such
  that $f(x,\cdot) \in \C_{\leqslant m} [T^*_xM]$ for any $x$.  
 If \eqref{eq:def_opdif_heis} holds on $U$, then $ \si (P)(x,\xi) =
 \sum_{|\al|\leqslant m } f_{0,\al}(x) \xi^{\sharp_x \al}$.
 \item[-] for any $P \in \dheis^m$,  $\si (P)
  = 0$ if and only if $P \in k^{-\frac{1}{2}} \dheis^{m-1}$. 
\item[-] for any $m \in \N$ and $p \in \R$, $\dheis^m \circ \pheis^p \subset
  \pheis^{m+p}$. Furthermore
  $$\si ( P \circ Q) (x, \cdot)  = \si (P)(x,\cdot) \, \sharp_x \, \si (Q) (x,
  \cdot)$$   for
  any $P \in \dheis^m$, $Q \in \pheis^p$. 
\end{itemize}
\end{prop}

\begin{proof}

{\cc   Recall that by Lemma  \ref{lem:corder}, if $P$ is a Heisenberg
  pseudodifferential operator of order $m$ on $U$ with principal symbol $\si$, then $\tilde
  \pi_j \circ P$ is a Heisenberg pseudodifferential operator of order $m+1$
  with principal symbol
  $$  (  \xi_j + \tfrac{i}{2}  \textstyle{\sum_{\ell=1}^n} \om_{j\ell } ( x)
  \partial_{\xi_\ell} ) \si (x,\xi) =  (\xi_j \sharp \si)(x,\xi)$$
  by \eqref{eq:EEEquation}. Then, starting from the fact that the identity is
  a Heisenberg pseudodifferential operator of order $0$ with principal symbol
  $1$, we deduce  by induction on $|\al|$ that $\tpi^{\al}$ is a
  Heisenberg pseudodifferential operator of order $|\al|$ with principal
  symbol $\xi^{\sharp \al}$. }


The first two
assertions follow. The third assertion is a consequence of  the fact that the $\xi^{\sharp
  \al} |_x$, $ \al \in \N^n$  are linearly independent so that $\sum f_{0, \al} \xi^{\sharp
  \al} = 0$ implies that $f_{0, \al} =0$. Last assertion follows again from
Lemma \ref{lem:corder} by induction on $m$. 
\end{proof}

\section{Resolvent} \label{sec:resolvent}

Let $(F, \la)$ be a real symplectic vector space with dimension $2d$. The Weyl product of the
Schwartz space $\mathcal{S} (F)$ is defined by
\begin{gather} \label{eq:prod_Weyl}
(a \circ_\la b ) (\xi) = (\pi)^{-2d } \int e^{-2i \la ( \eta , \zeta) }
a(\xi+\eta) b ( \xi + \zeta) \; d\mu_F( \eta) \;
d\mu_F ( \zeta)
\end{gather}
where $\mu_F = \la^{\wedge d} /d! $ is the Liouville measure  of $F$.  
For $F = \R^d \times \R^d$ with
\begin{gather} \label{eq:la-t-tau}
\la (t,\tau; s,  \varsi) = \tau \cdot s  -
\varsi \cdot t ,  \qquad ( t, \tau), (s, \varsi) \in \R^{d} \times \R^d
\end{gather}
\eqref{eq:prod_Weyl} is the
composition law of the Weyl symbols of pseudodifferential operators of $\R^d$,
cf. for instance \cite[page 152]{Ho}.

This product extends continuously from $S^m (F) \times S^p (F)$ to $S^{m+p} (
F)$ by preserving the subspace of polyhomogeneous symbols.  So the
corresponding pseudodifferential operators, $f^w ( x, \frac{1}{i} \partial)$,
with $f \in S^\infty ( \R^{2d})$,  form an algebra, called sometimes 
the Shubin class or isotropic {\cc algebra}. This algebra is one of the most studied
in microlocal analysis,  cf.
\cite[Chapter IV]{Sh}, \cite{He}, \cite[Chapter 4]{Me}, \cite[Chapter 4]{EM},
\cite[Appendix A]{Ta} for lecture note references.

The Weyl product  appears naturally in our context as the product of the
operators $f( \frac{1}{i} \nabla^A )$ defined in Section
\ref{sec:isotropic-algebra} when $A$ is nondegenerate. 
\begin{lemme} \label{lem:Weyl_algebra}
  If $A \in \wedge^2 E^*$ is nondegenerate, then for any $f$, $g$
  in $S^{\infty} ( E^*)$, $$ f( \tfrac{1}{i} \nabla^A) \circ g ( \tfrac{1}{i}
  \nabla^A) = ( f \circ_{\la} g ) ( \tfrac{1}{i} \nabla^A)$$ where $\la$ is the
  symplectic form of $E^*$ dual to $A$.
  \end{lemme}
  \begin{proof}
Introduce a symplectic basis $(e_i, f_i)$ of $(E, A)$ and denote by $(x_i, y_i)$ 
the asso\-cia\-ted linear coordinates, so that $E = \R^d_x \times \R^d_y$. Then the operators
$$\tfrac{1}{i}
\nabla_{e_i}  = \tfrac{1}{i} \partial_{x_i} + \tfrac{1}{2} y_i, \quad
\tfrac{1}{i} \nabla_{f_i} =  \tfrac{1}{i} \partial_{y_i} - \tfrac{1}{2} x_i,
\quad \  \tfrac{1}{i} \partial_{y_i} +
\tfrac{1}{2} x_i, \quad    \tfrac{1}{i} \partial_{x_i} -
\tfrac{1}{2} y_i $$
satisfy the same
commutation relations as the operators $ s_i $, $ \frac{1}{i} \partial_{s_i}$,
$t_i$, $\frac{1}{i} \partial_{t_i}$ of $\R^d_s \times \R^d_t$. So the linear
isomorphism $\Phi : \R^{4d} \rightarrow \R^{4d}$,
$$ \Phi ( x, \xi, y, \eta ) = ( \xi + \tfrac{1}{2} y,\, \eta - \tfrac{1}{2} x,\,
\eta + \tfrac{1}{2} x,\, \xi - \tfrac{1}{2} y ) .$$
is a symplectomorphism. The  metaplectic representation yields us a unitary
operator $U : L^2 (E)
\rightarrow L^2 ( \R^{2d})$ satisfying
$$ f^w  = U (f \circ \Phi ) ^w U^*, \qquad \forall \, f \in \mathcal{S}' ( \R^{4d}),$$   cf. \cite[Theorem 18.5.9]{Ho}. Applying this to $f (x, \xi, y, \eta ) = g ( \xi + \frac{1}{2} y, \eta - \frac{1}{2} x
) =  ( ( g \boxtimes 1 ) \circ \Phi ) ( x, \xi, y, \eta)$, we obtain 
$$  f ^w ( x, \tfrac{1}{i} \partial_x, y, \tfrac{1}{i}
\partial_y)  = U ( g^w ( s, \tfrac{1}{i} \partial_s) \otimes \op{id}_{\R^d_t})U^* $$
and by Lemma \ref{lem:Weyl_symbol}, $ f ^w ( x, \tfrac{1}{i} \partial_x, y, \tfrac{1}{i}
\partial_y) =  g ( \tfrac{1}{i} \nabla^A ) $. The result follows. 
\end{proof}

From now on, we assume that $\om$ is nondegenerate. By Lemma
\ref{lem:Weyl_algebra}, at any $x \in
M$, the product $\sharp_x$ defined in  \eqref{eq:sharp_produit} extends continuously 
$$S^{m} ( T_x^*M)
\times S^p  (T_x^*M) \rightarrow S^{m+p}  (T_x^*M). $$
{\cd{We are now ready to consider the spaces $S^{m} ( M,T^*M)$ of symbols defined
on $T^*M$, cf. \eqref{eq:symb_fibre}.}} We say that $f \in S^{m} (M, T^*M)$ is  
elliptic if $|f(x,\xi)|\geqslant C ^{-1} |\xi|^m$ when $| \xi|
\geqslant C$ for some positive $C$. We say that $f$ is invertible if at any $x \in M$, $f(x,
\cdot)$ is invertible in $(S^\infty (T_x^*M) , \sharp_x)$. 

\begin{lemme} \label{lem:prod_inverse_symbol}
  \begin{enumerate}
    \item
$S_{\op{ph}}^\infty ( M, T^*M)$ endowed with the fibered product $( f \sharp g ) (x) = f
( x) \sharp_x g ( x) $ is a filtered algebra.
\item For any $ f \in S^m
_{\op{ph}} (M, T^*M)$ which is both elliptic and invertible, the pointwise
inverse of $f$ belongs to $S^{-m}
_{\op{ph}} ( M, T^*M)$. 
\end{enumerate}
%
%
\end{lemme}

\begin{proof} This holds more generally for $S^m_{\op{ph}} ( N,E)$ where $E$ is
  any symplectic vector bundle with base $N$. When $N$ is a point, $S^\infty_{\op{ph}} (N,E)$ 
  is isomorphic with the Weyl algebra $S^\infty_{\op{ph}} ( \R^{2d})$, and the result is
  well-known as we already mentioned it.

  In general, we can assume that $E$ is the trivial symplectic bundle $\R^{2d}$ over
  an open subset $U$ of an Euclidean space, {\cc so that the product $\sharp_x$ is independent
    of $x \in U$. So the first assertion is that the Weyl product  \eqref{eq:prod_Weyl} with $F = \R^{2d}$,
$\la$ given by \eqref{eq:la-t-tau} and symbols $a$ and $b$ depending
smoothly on an additional parameter $x \in U$, is continuous 
\begin{gather} \label{eq:continuite} 
      S^{m} ( U, \R^{2d})  \times S^{p} (U, \R^{2d}) \rightarrow S^{m+p} ( U,
      \R^{2d}).
    \end{gather}
We claim that this follows from the particular case where $U$ is reduced to a
point, which is well-known as already mentioned. To prove the claim, we will use
 the following easy facts:  if $f \in S^m ( U, \R^n)$, then $x \rightarrow f (x, \cdot)$ is
 continuous from $U$ to $S^m ( \R^n)$. Conversely, if $x \rightarrow f(x,
 \cdot )$ is continuous from $U$ to $S^m ( \R_\xi^n)$, then the partial
 derivatives $\partial_\xi^\al f (x, \xi)$ depend continuously on $(x, \xi) \in U
 \times \R^n$. Now
 consider any bilinear continuous map
 $ B : S^m (  \R^n) \times S^p (  \R^n) \rightarrow S^{\ell} ( \R^n)$. Let us
 show that the associated map
\begin{gather} \label{eq:Betendue}
 S^m (U,  \R^n) \times S^p ( U , \R^n) \rightarrow S^{\ell} (U, \R^n) 
\end{gather}
sending $(f,g)$ to $h$ given by $h (x, \xi ) = B( f(x, \cdot), g(x, \cdot) ) (\xi)$ is
 well-defined and continuous. By the preliminary observation, the functions
 $\partial^\al_\xi h$ are all continuous. Let us prove that $h$
 is derivable with respect to $x$ as well. By Taylor expanding we have
 $$ f( x + u, \xi ) = f (x, \xi ) + \textstyle{\sum}u_i \partial_{x_i} f ( x,
 \xi) + \psi (x, \xi, u) $$
and the usual integral formula of the remainder shows that $\psi (\cdot, u ) \in S^m ( U, \R^n)$ for any $u$, with all its
 semi-norms in $\bigo ( |u|^2)$. Expanding similarly $g$, we get that
  $h( x + u , \cdot )$ is equal to $$  h (x, \cdot) + \textstyle{\sum} u_i \bigl( B (
 \partial_{x_i} f (x, \cdot), g( x, \cdot) ) + B ( f ( x, \cdot) ,
 \partial_{x_i} g ( x, \cdot)) \bigr) + \bigo ( |u|^2)$$ 
where the $\bigo$ is for all the seminorms of $S^{\ell} ( \R^n)$. It follows that $h$ is derivable with respect to $x$, with continuous partial
 derivatives given by
 $$ \partial_{x_i} h ( x, \xi ) =  B (
 \partial_{x_i} f (x, \cdot), g( x, \cdot) ) (\xi)  + B ( f ( x, \cdot) ,
 \partial_{x_i} g ( x, \cdot))(\xi) . $$
 Since $\partial_x^\al f \in S^m ( U, \R^n)$ and similarly for $g$, it follows
 by induction that $h$ is smooth and
 \begin{gather} \label{eq:partial_deriv} \partial_{x}^\ga h ( x, \xi ) = \textstyle{\sum}_{\al + \be = \ga} C_{\al, \be} B (
 \partial_x^\al f ( x, \cdot) , \partial_x^\be g ( x, \cdot ) ) (\xi) 
\end{gather}
So $h$ belongs to $S^\ell ( U, \R^n)$. The continuity of \eqref{eq:Betendue} follows easily from \eqref{eq:partial_deriv}.


  }


    Let  $f \in
S_{\op{ph}}^{m}(U, \R^{2d})$ be elliptic and
invertible. Let us prove that its pointwise inverse $g$ is in $S^{-m}
_{\op{ph}} ( U, \R^{2d} )$. Multiplying $f$ by
$f(x_0)^{-1}$, we may assume that $m =0$. Since $S^\infty_{\op{ph}} ( U, \R^{2d})$ is a
filtered algebra, cf. \eqref{eq:continuite},  and by
Borel lemma, $f$ has a parametrix $h \in S^0 _{\op {ph}} (U, \R^{2d})$. Thus
$$ h \sharp f = 1 + r, \qquad f \sharp h = 1 + s  \qquad \text{with }  r, \, s
\in S^{-\infty} ( U, \R^{2d}).$$
Let us
prove that $g = h + S^{-\infty} ( U, \R^{2d})$.
 We have 
$g  =
      h - r \sharp h + r \sharp g \sharp s$. By 
      \eqref{eq:continuite} again, $r \sharp h \in  S^{-\infty} ( U,
      \R^{2d})$. It remains to prove that $r \sharp g \sharp s \in  S^{-\infty} ( U,
\R^{2d})$.

By Calderon-Vaillancourt theorem, the Weyl quantization $\op{Op}: S^{0} (
\R^{2d} ) \rightarrow \mathcal{L} (L^{2 } ( \R^d))$ is continuous. $\op{Op}
(g (x) )$ being the inverse of $\op{Op} ( f(x))$ for any $x$,  $\op{Op}
(g) \in \Ci ( U , \mathcal{L} ( L^{2 } ( \R^d)))$.
We claim that  the multilinear map
\begin{gather} \label{eq:multilinear}
M: S^{-\infty} ( \R^{2d})  \times \mathcal{L} (
      L^2 ( \R^d) )\times S^{-\infty} ( \R^{2d})  \rightarrow S^{-\infty} (
      \R^{2d}),
    \end{gather}
    defined by $\op{Op} (
      M(\si ,A, \tau)) = \op{Op} ( \si
      )\circ A \circ \op{Op} (\tau
      ) $, is continuous, which implies that that   $r \sharp g \sharp s = M ( r, \op{Op}( g) ,  s)$ belongs to  $S^{-\infty} ( U,
      \R^{2d})$. {\cc To prove the claim, recall first that Weyl quantization is an
      isomorphism between $S^{-\infty} (\R^{2d}) = \mathcal{S} ( \R^{2d}) $ and the space of linear maps having a
      Schwartz kernel in $\mathcal{S} ( \R^{2d})$, that is the space of linear
      continuous maps $\mathcal{S}' ( \R^d) \rightarrow  \mathcal{S} (\R^d)$.
      So for any  $\si$, $\tau \in \mathcal{S} ( \R^{2d})$ and
      $A$ linear continuous $\mathcal{S} (\R^d) \rightarrow \mathcal{S}' (
      \R^d)$, the composition $\op{Op} ( \si
      )\circ A \circ \op{Op} (\tau
      ) $ is well-defined and has the form $\op{Op} ( \rho)$ with $\rho \in
      \mathcal{S} ( \R^{2d})$. Moreover, the continuity of \eqref{eq:multilinear} is
      equivalent to the continuity of
      $$ M' : \mathcal{S} ( \R^{2d} ) \times \mathcal{L} (
      L^2 ( \R^d) )\times \mathcal{S} ( \R^{2d} ) \rightarrow \mathcal{S} (
      \R^{2d} ) $$
      defined by $M' (S,A, T) = S\circ A \circ T$ where we identify Schwartz
      kernels with their associated linear map. An explicit formula  in
      terms of scalar product of $L^2( \R^d)$ is: 
 \begin{gather} \label{eq:mprimeformule}
       M' ( S,A, T) (x,y) = ( A T( \cdot, y ) , \con{S  (x, \cdot )} )_{L^2
         ( \R^d)} . 
       \end{gather}  
      Introduce the norm $\| S \|_m = \sup_{(x,y) \in \R^{2d}} |S(x,y)|
      \langle x,y \rangle ^m  $ for positive $m$. Since $\langle x \rangle
      \leqslant \langle x,y \rangle$ and $
      \langle y \rangle \leqslant \langle x,y \rangle $,  we have
 \begin{gather} \label{eq:forint}
       |S(x,z)| \langle x
       \rangle ^m \leqslant \langle z \rangle^{-p} \| S \|_{m+ p }, \qquad |
       T(z,y) | \langle y \rangle^m \leqslant \langle z \rangle ^{-p} \| T\|_{m+p} .
     \end{gather}
    Choose $p > n/2$ so that $\langle \cdot \rangle ^{-p}$ is in $L^2 (
      \R^d)$. Then using that $\langle x, y \rangle \leqslant \langle x
      \rangle \langle y \rangle$, it follows from \eqref{eq:mprimeformule} and
      \eqref{eq:forint} that  
$$    \| M' ( S,A, T) \|_{m+p}  \leqslant C  \| A \| \, \|
S \|_{m + p } \| T \|_{m + p }$$
The estimates of the derivatives are similar.}     \end{proof}

Consider now $P \in \dheis^m ( L, \nabla)$ having an elliptic symbol. Then for any
fixed $k$, $P_k$ is an elliptic differential operator of $\Ci ( M, L^k)$, so
for any $s \in \R$, $P_k$ extends to a Fredholm operator of $\mathcal{L} (
H^s (M, L^k) , H^{s -m } ( M , L^k))$.  If we assume that the symbol of $P$ is
invertible, then by the following Theorem, $P_k$ is invertible when $k$ is
large, and its inverse is a Heisenberg pseudodifferential operator.

  \begin{theo} \label{theo:resolvent}
Assume that $\om $ is nondegenerate.  Let $P \in \dheis^m (L, \nabla)$ having an elliptic and invertible symbol $\si \in \Ci ( T^*M)$. Then
    there exists $ Q \in \pheis^{-m} (L, \nabla)$ such that
    \begin{enumerate}
      \item[-] $PQ -\op{id} $ and $QP -
        \op{id}$ are in $k^{-\infty} \Psi^{-\infty} (L)$
        \item[-]  when $k$ is
    sufficiently large, $ Q_k P_k = P_k Q_k =
    \op{id}$
    \item[-] the symbol of $Q$ is  the
        inverse    of $\si$ for the product $\sharp$.
    \end{enumerate}
  \end{theo}

\begin{proof} This follows merely from the previous results, by the standard
  techniques for elliptic operators.  First, using Lemma
  \ref{lem:prod_inverse_symbol}, we construct  a
      parame\-trix $Q \in \pheis^{-m} (L, \nabla)$ of $P$, so  $PQ = \op{id} + R  $ and $QP = \op{id} + S $ with $R, S$ in
      the residual algebra $k^{-\infty} \Psi ^{-\infty} ( L)$.  Then, by the
      Sobolev continuity \eqref{eq:soblov_semi_classical}, $R_k$ and $S_k$
      belongs to $\mathcal{L} ( L^2 ( M, L^k))$ and their operator norms are in $\bigo
      ( k^{-\infty})$.  So when $k
      \geqslant k_0$, $P_k$ is invertible from $H^{m} ( M , L^k)$ to $H^0
      (M, L^k)$, which implies by the Fredholm properties of elliptic
      operators \cite[Theorem 8.1]{Sh}, that  $P_k$ is an invertible
      operator of the distribution space $\mathcal{D}' ( M , L^k)$. 

   Its inverse satisfies 
\begin{gather} \label{eq:eqint}
      P^{-1} _k = Q_k - R_k Q_k + R_k (Q_k P_k)^{-1} Q_k S_k
    \end{gather}
    By Lemma \ref{lem:residuel}, $(R_k Q_k)$ and $(Q_k S_k)$ are in $k^{-\infty}
\Psi^{-\infty} (L)$. It is a classical fact  that if $(A_k)$, $(B_k)$ are in
$k^{-\infty} \Psi^{-\infty} (L)$ and $C_k = \bigo (1): L^2 (M, L^k) \rightarrow
L^2 ( M , L^k)$, then $(A_k C_k B_k)$ is in $k^{-\infty} \Psi ^{-\infty} (L)$.
So the last term in \eqref{eq:eqint} belongs to $k^{-\infty} \Psi^{-\infty}
(L)$. So by adding to $Q_k$ an element of the residual algebra, we have that $Q_k = P_k^{-1}$ when $k$ is large. 
\end{proof}

Assume $m>0$ and consider $P \in \dheis^m ( L, \nabla)$ having an elliptic symbol $\si$ such that for
some $z_0 \in \C$, $\si - z_0$ is invertible. Then by
Theorem \ref{theo:resolvent}, when $k$ is sufficiently large, $P_k -z_0 $ has an
inverse, which is continuous $ L^2 (M, L^k) \rightarrow H^m ( M, L^k) $. So
the restriction of $P_k$ to $H^m (M, L^k)$ is a closed unbounded operator of $L^2 ( M, L^k)$ having a compact resolvent. So its spectrum is a discrete subset of
$\C$ and it consists only of eigenvalues with finite multiplicity, the
generalised eigenvectors being smooth \cite[Theorem 8.4]{Sh}.

To state the next theorem, we need some spectral properties of
the symbols themselves. Later we will explain these properties in
terms of Weyl quantization, but since this quantization is only auxiliary in
what we do, we prefer first to discuss everything intrinsically in terms of
the algebra
$(S^\infty(F), \circ_{\la} )$ where $(F, \la)$ is a symplectic vector space as
above. 

The spectrum of $a \in S^\infty(F)$ is defined by:  $z \notin \op{sp} ( a) $ if
and only if $z- a$ is invertible in $(S^{\infty} (F), \circ_\la )$. A family $(b(z), z \in \Om)$ of $S^m (F)$ is
holomorphic if $\Om$ is an open set of $\C$,  $b \in S^m (\Om, F)$
and $ \partial_{\con{z}} b =0$.
By the analytic Fredholm
theory {\cc for the operators with symbols in $S^{\infty} (F)$} exposed in \cite[Chapter 3]{Me}, 
 for any elliptic $a \in
S^m_{\op{ph}} (F)$ with $m>0$, the spectrum of $a$ is $\C$ or a discrete
subset of $\C$. In the latter case, the resolvent $ ( (a - z)
^{(-1)_{\circ_\la}}, \; z \in \C \setminus \op{sp} (a))$
is a holomorphic family of $S^{-m} ( F)$ and for any $z_0 \in \op{sp} (a
)$, we have on a 
neighborhood of $z_0$ for some $N \in \N$
\begin{gather}  \label{eq:residu}
(a -z) ^{(-1)_{\circ_\la}} = h(z) + \frac{r_1}{z-z_0} +
\ldots + \frac{r_N}{(z- z_0)^N}
\end{gather}
where $(h(z))$ is a holomorphic family of $S_{\op{ph}}
^{-m} (F)$   and $r_1$, \ldots, $r_N$ are
in $S^{-\infty}(F)$.

\begin{theo} \label{theo:resolvent_et_projecteur}
  Assume that $\om$ is nondegenerate. Let $P \in \dheis^m ( L ,
  \nabla)$  be elliptic with $m \geqslant 1$ and symbol $\si$. 
  Let  $\Sigma$ be the closed set $\bigcup_{x
    \in M } \op{sp} (\si (x))$. Then
  \begin{enumerate}
    \item   if $K$ is a compact subset of $\C$ disjoint from $\Sigma$, then the spectrum of $P_k$
  does not intersect $K$ when $k$ is large enough.
\item if $\Om $ is an open bounded subset of $\C$  with a smooth boundary disjoint
  from $\Si$, then there exists $\Pi \in \pheis^{-m}
  (L, \nabla)$ such that $\Pi_k =   (1_{\Om} ( P_k))$ when $k$ is large. Furthermore the principal symbol of $\Pi$ is at $x$
\begin{gather} \label{eq:symbolpi}
  \pi
    (x) = (2i\pi)^{-1} \int_{\partial \Om}  (\si (x)- z)^{(-1)_{\sharp_x}} dz. 
  \end{gather}
  \item if for any $k$, $P_k$ is formally self-adjoint for some volume element
    of $M$, then for any $E_-, E_+
    \in \R \setminus \Si$ with $E_- < E_+$, $(1_{[E_-,
      E_+]} ( P_k) )  $ belongs to $\pheis^{-\infty}  ( L , \nabla)$. 
\end{enumerate}
\end{theo}

Observe that the symbol $\pi(x)$ is
the sum of the residues of the poles in $\Om$  of the resolvent of $\si (x)$.
As we will see in the proof, the third assertion is a particular case of the
second one, the symbol being the sum of the residues of the poles in
$[E_-,E_+]$. 

In \cite{oim_Heis}, we will prove that $\pheis^m \circ \pheis^p \subset
\pheis^{m+p}$. So in the second assertion, $\Pi$ being idempotent, it belongs to $\pheis^{-\infty} (L,
\nabla)$.


\begin{proof} 
First, $\Si $ is closed because the Weyl quantization is continuous from
$ S^m (\R^d)$ to $ \mathcal{L} ( H^m_{\op{iso}} (\R^d), H^0 _{\op{iso}} (
\R^d) )$, so that the characterization of the spectrum given below implies
that if $z_0 \notin \op{sp} ( \si(x_0))$ then $z \notin \op{sp} ( \si 
(x))$ when $(z,x)$ is sufficiently close to $(z_0 , x_0)$.

Assume that $K$ is a compact subset of $\C$ disjoint from $\Si$.  When $z \in
K$,  $(P_k - z) $ satisfies the assumptions of Theorem \ref{theo:resolvent},
so there exists $Q (z) \in \pheis^{-m} (L, \nabla) $ such
that $Q_k (z) = (z - P_k)^{-1}$ when $k \geqslant k_0 (z)$.
This proves at least that $P_k$ has a
compact resolvent as explained above when $k$ is large. Moreover we claim that everything in the proof
of Theorem \ref{theo:resolvent} can be done continuously  with respect to $z \in K$ (even holomorphically
with respect to $z$ in a neighborhood of $K$). More precisely, the Schwartz
kernel of $Q_k (z) $ is locally of the form \eqref{eq:Weyl_form} where
the dependence in $z$ in only in the symbol $b$, which is continuous in $z$. This proves first that we
can choose $k_0 (z) $ independent of $z$, which shows the first assertion.
Second, if $\Om$ satisfies the assumptions of the second assertion of the Theorem, we can
apply the previous consideration to $K = \partial \Om$ and it follows that $  \Pi_k :=  (2 i  \pi)^{-1} \int_{\partial \Omega } Q_k (z) \; dz $
belongs to $\pheis^{-m} (L, \nabla)$ with a symbol given by
\eqref{eq:symbolpi}. When $k$ is large enough, $Q_k(z)$ is the resolvent, so
by Cauchy formula, $\Pi_k = 1_{\Om} ( P_k)$. This concludes the proof of the second assertion.

For the last assertion, by assumption, for any fixed $k$, $P_k$ is a formally self-adjoint
elliptic differential operator on a compact manifold, so its spectrum is a discrete
subset of $\R$ and $1_{[E_-, E_+]} (P_k)$ is a finite rank projector onto a
subspace of $\Ci ( M , L^k)$, \cite[Theorem 8.3]{Sh}.  Moreover $\si$ is real valued so $\Si
\subset \R$. So there exists $\Om$ satisfying the previous assumptions and
such that $\Om \cap \R = [E_-, E_+]$. So $ \Pi = (1_{[E_-, E_+]} (P_k) , \; k \in \N)$ belongs to
$\pheis^{-m} ( L, \nabla)$.

For any odd $N \in \N$, $\Pi = 1_{[E^N_-, E^N_+]}
(P^N)$ and $P^N \in \dheis^{mN} ( L , \nabla)$,  which implies by the previous
argument that $\Pi$ belongs to $ \pheis^{- N m} (L, \nabla) $, so $\Pi \in \pheis^{- \infty}(L, \nabla)$. 
\end{proof}

Let us discuss briefly the invertibility and resolvent of elliptic elements of
$(S^\infty (F), \circ_\la) $ from the point of view of Weyl quantization.
Let $\piso^{\infty} ( \R^d)$ be the space of pseudodifferential operators
of $\R^d$ with a symbol in $S^\infty_{\op{ph}} (\R^{2d})$. Any $A \in \piso^m
( \R^d)$ acts continuously  $\mathcal{S}  ( \R^d)\rightarrow \mathcal{S}  ( \R^d)$,  
$\mathcal{S}' ( \R^d) \rightarrow \mathcal{S}' ( \R^d)  $ and $H^s_{\op{iso}}  ( \R^d)
\rightarrow H^{s-m}_{\op{iso}}  ( \R^d)$, where
$H^s_{\op{iso}}  ( \R^d)$ are the
isotropic Sobolev spaces  
 $$H^s_{\op{iso} }  ( \R^d) = \{ u \in \mathcal{S}' ( \R^d) , \;
 A u \in L^2 ( \R^d), \; \forall A \in \piso^s ( \R^d) \}, \quad s \in \R.$$
 When $A$ is elliptic, the following Fredholm property holds:
$\ker A$ and $\ker A^*$ are finite dimensional subspaces of $\mathcal{S}  ( \R^d) $,
$$ \mathcal{S}' ( \R^d)  = A( \mathcal{S}' ( \R^d)) \oplus \ker A^* = A^* ( \mathcal{S}' ( \R^d)) \oplus
\ker A$$ and the generalised inverse $B : \mathcal{S}' ( \R^d) \rightarrow
\mathcal{S}' ( \R^d)$ such that $BA - \op{id}$ and $AB - \op{id}$ are the orthogonal
projectors onto $\ker A$ and $\ker A^*$ respectively, belongs to $\piso^{-m} (
\R^d)$. So $A$ is invertible in the algebra $\piso^{\infty} ( \R^d)$ if and
only if 
$\ker A = \ker A^* =0$ if and only if  $A$
is invertible as an operator in $\mathcal{S}'$ if and only if  $A$ is invertible in
$\mathcal{L} ( H^s_{\op{iso}}  ( \R^d), H^{s-m}_{\op{iso}}  ( \R^d) )$.

When $m>0$, any elliptic $A \in \piso^m ( \R^d)$ defines by restriction a
closed unbounded
operator of $L^2 ( \R^d)$ with domain $H^m_{\op{iso}}  ( \R^d)$. By the previous
characterization of invertibility, the spectrum of $A$ is the same as the 
spectrum of its symbol $a$ defined above. Assume it is not empty, then
$A$ has a compact resolvent, and as it was already explained, 
$\op{sp} (A)$ is a discrete subset of $\C$ and the resolvent $(A-z)^{-1}$ is a
holomorphic family of $\piso^{-m} ( \R^d)$.
Furthermore, for $A = a^w (x, \frac{1}{i} \partial_x )$,  the residues $r_{\ell}^w
(x, \frac{1}{i} \partial_x)$ defined in \eqref{eq:residu} have finite rank and $r_{1}^w
(x, \frac{1}{i} \partial_x)$ is a projector onto the space of generalised
eigenvectors of $A$ for the eigenvalue $z_0$, which is a subspace of
$\mathcal{S}  ( \R^d)$.

{\cc To end this section let us prove a Weyl law corresponding to Theorem
\ref{theo:resolvent_et_projecteur}. {\cd{Consider again a real symplectic vector
space $( F, \lambda)$ and the associated algebra $(S^{\infty} (F),
\circ_{\lambda})$}}. Then the Schwartz class $
\mathcal{S} (F)= S^{-\infty} (F)
$ is an ideal of $(S^{\infty} (F), \circ_{\lambda})$. Set
\begin{gather} \label{eq:def_trace}
  \op{tr} a := \int_F a \, d\mu_{\la}, \qquad \forall a \in \mathcal{S} (F) 
\end{gather}
where $\mu_{\la}$
is the Liouville measure of $(F, \lambda)$. Then by \cite[Chapter 3.14]{Me}, $\op{tr}$ is a trace in the
sense that $\op{tr} (ab) = \op{tr} ( ba)$ for any $a \in \mathcal{S}(F)$ and
$b \in S^{\infty} (F)$. Moreover, for $F = \R^{2d}_{x,\xi}$
with $\lambda$ the
usual symplectic form, by \cite[Section 27.1]{Sh}, for any $a \in \mathcal{S} (F)$, $a^{w} (
x, \frac{1}{i} \partial_x )$ is a trace class
operator of $L^2 ( \R^d)$ and $\op{tr} a^{w} ( x, \frac{1}{i} \partial_x ) =
\op{tr} a$. In particular, when $  a^{w} ( x, \frac{1}{i} \partial_x )$ is a
projector, it has a finite rank equal to $\op{tr} (a)$.

Consider now $Q \in \Psi ^{-\infty}_{\op{Heis}} ( L, \nabla )$. Then for any
$k$, $Q_k$ is a smoothing operator, so  $Q_k$  is trace class operator of $L^2
( M, L^k)$ and its trace is the integral of the Schwartz kernel on the diagonal.
Observe that the Schwartz kernel being a section of $(L^k\boxtimes \con{L}^k)
\otimes ( \C_M \boxtimes |\Lambda| (M) )$, its restriction to the diagonal
identifies naturally with a section of $|\Lambda|( M)$, so its integral is
well-defined. By the expression \eqref{eq:def_heis_global} for the Schwartz
kernel of a Heisenberg pseudodifferential operator, it comes that
\begin{gather} \label{eq:trace}
\op{tr} (Q_k) = \Bigl( \frac{k}{2 \pi} \Bigr)^{d} \int_{T^*M} \si (Q) \,
d\mu_{T^*M} + \bigo ( k^{d-1/2} )
\end{gather}
where $\mu_{T^*M}$ is the Liouville measure of $T^*M$. Formula
\eqref{eq:trace} holds without assuming that the curvature of $\nabla$ is
non-degenerate, so $d = n/2$ is not necessarily an integer. If $\om$ is symplectic, then we have
\begin{gather} \label{eq:int_trace}
  \int_{T^*M} \si (Q) \,
  d\mu_{T^*M} = \int_M \op{tr} \si(Q) (x, \cdot ) \; d \mu_M (x)
\end{gather}
where $\op{tr} \si(Q) (x)$ is the trace \eqref{eq:def_trace} of $\si(Q) (x) \in \mathcal{S} (T_x^*M)$ and $\mu_M$ is the Liouville measure of $(M,
\om)$. The proof of \eqref{eq:int_trace} is straightforward in Darboux
coordinates of $(M,\om)$. Assume now that for any $x \in M$, $\si (Q) (x )$
is a projector. Then its trace is an integer depending continuously of $x$, so
it is constant equal to $N \in \N$ and
$$ \int_M \op{tr} \si(Q) (x) \; d \mu_M (x) = N \op{Vol} ( M, \mu_M) $$    
Applying this to $Q_k =  (1_{[E_-,E_+]} ( P_k) )$ defined as in Theorem
\ref{theo:resolvent_et_projecteur}, we obtain the following Weyl estimate.

\begin{cor} \label{cor:Weyl_law} For any $P$, $E_-$, $E_+$ satisfying the assumptions of the third assertion of Theorem
  \ref{theo:resolvent_et_projecteur},  we have
  $$ \sharp \bigl( \op{sp}(P_k) \cap [E_-, E_+ ] \bigr) = \Bigl( \frac{k}{2 \pi}
  \Bigr)^{d} N \op{ Vol} ( M , \mu_M) + \bigo ( k^{d-\frac{1}{2}} )$$
  where $N = \sharp \bigl( \op{sp} (\si(P) (x))  \cap [E_-, E_+ ] \bigr) $ for any $x \in M$.
\end{cor}
}

\section{Auxiliary bundles} \label{sec:auxiliary-bundles}

Let us first define symbols taking values in an auxiliary bundle.
Recall the spaces $S^m_* (N, E)$ introduced in Section
\ref{sec:twist-pseud-oper} for a real vector bundle $p:E \rightarrow N$ and $*
= \emptyset$, $\op{ph}$, $\op{sc}$. Let $B $ be a complex vector bundle over
$N$. By definition $S^m  ( N,E; B)$ is the space of
sections $s \in \Ci ( E,
p^*B)$ such that for any frame $(u_{\al})$ of $ B$ over an open set $U$ of $N$,
we have over $p^{-1} (U)$,
$$ s( x, \xi ) = \sum f_{\al}  ( x, \xi ) u_{\al} (x) , \qquad x \in N, \; \xi
\in E_x $$
with coefficients $f_{\al}$ in $S^m ( U, E)$. Since $S^m (U, E)$ is a $\Ci (
U)$-submodule of $\Ci ( U, E)$, this definition is compatible with the frame
changes. Similarly, we define $S^m_{*} (M,E;B)$ for $*  = \op{ph}$ or
$\op{sc}$  by requiring that the
coefficients $f_{\al}$ belong to  $S^m_{*} (U,E)$. More precisely, in the case
of semiclassical symbols where the section $s$ and its local coefficients depend on
$h$, we only choose frames $(u_{\al})$ independent of
$h$.

Let $A_1$ and $A_2$ be two complex vector bundles over $M$ and let us define
the pseudodifferential operator spaces $\Psi^m_{\op{sc}} (M; A_1, A_2)$, $\ptw^m
(L; A_1, A_2)$ and $\pheis^m (L, \nabla; A_1, A_2)$. For $A_1$,
$A_2$ being both the trivial line bundle, these are the spaces we introduced
previously. In general, set $B = A_2 \boxtimes A_1^*$. Then  
\begin{itemize}
  \item  $\Psi^m_{\op{sc}} (M; A_1, A_2)$ consists of the families $(P_h : \Ci ( M , A_1)
\rightarrow \Ci ( M , A_2)$, $h \in (0,1])$ satisfying the same conditions as
before except that the amplitude $a$ appearing in
\eqref{eq:def_pseudo_semiclassique} belongs to $S_{\op{sc}} ^m (U^2, \R^n; B)$. 
\item  $\ptw^m ( L;A_1, A_2)$ consists of the  families
  \begin{gather} \label{eq:familyP}
    P= (P_k : \Ci ( M , L^k \otimes A_1) \rightarrow \Ci ( M ,L^k \otimes
    A_2),\; k \in \N) 
  \end{gather}
  satisfying the conditions of Definition \ref{def:twisted_pseudo} with   $a  \in S_{\op{sc}} ^m (U^2, \R^n; B)$
\item  $\pheis^m ( L, \nabla ;A_1, A_2)$ consists of the families $P$ of the form
  \eqref{eq:familyP}  satisfying the conditions of Definition \ref{def:heis-semicl-oper} with
  $Q_h$ an operator of $S^m_{\op{sc}} ( M; A_1, A_2)$
\end{itemize}
The symbol of $P$ is defined as before. Since the restriction of $B$ to the
diagonal is isomorphic with $\op{Hom} (A_1, A_2)$, in the three cases, the
symbol identifies with an element of 
$S^m_{\op{ph}} ( M, T^*M; \op{Hom} ( A_1, A_2))$.

The space $\dheis^m ( L,\nabla ; A_1, A_2)$ of Heisenberg differential
operators consists of
the families  \eqref{eq:familyP} of
differential operators such that for
any coordinate chart $(U, x_i)$ of $M$, we
have on $U$
\begin{gather} \label{eq:def_opdif_heis_aux} 
 P_k = \sum_{\ell \in \N, \; \al \in \N^n, \; \ell + |\al| \leqslant m} k^{-
  \frac{\ell}{2}} f_{\ell, \al } \tpi^\al
\end{gather}
where $f_{\ell, \al} \in \Ci (U, \op{Hom}(A_1, A_2))$, $\tpi^{\al} =
\tpi_1^{\al (1)} \ldots \tpi_n^{\al (n)}$ with $\tpi_i=\tfrac{1}{i \sqrt k}
\nabla_{\partial_{x_i}} ^{L^k\otimes A_2}$. Here we use a connection of $A_2$,
which induces with the connection of $L$ a covariant derivative of $A_2
\otimes L^k$. Proposition
\ref{prop:heis-diff-pseudo} still holds:  Heisenberg differential operators
are Heisenberg pseudodifferential operators, the symbol of
\eqref{eq:def_opdif_heis_aux} is $ \sum_{|\al|\leqslant m}  f_{0,\al} ( x)
\xi^{\sharp_x \al}$, 
$$ \dheis^m ( L, \nabla; A_2, A_3) \circ \pheis^p ( L, \nabla; A_1, A_2)
\subset \pheis ^{m+p} (L, \nabla ; A_1, A_3),$$
and the product of symbols  is the fiberwise product $\sharp_x$ tensored by the composition
$\op{Hom} ( A_{2,x}, A_{3,x}) \times \op{Hom} ( A_{1,x}, A_{2,x})  \rightarrow \op{Hom} ( A_{1,x}, A_{3,x}) . $
It is easy to see  that the definition of the Heisenberg differential operators and of their
symbols do not depend on the choice of the connection  of $A_2$.

In the sequel we assume that $A_1 =A_2 =A$ and is equipped with a Hermitian metric. We use the notation $\dheis^m
( L, \nabla; A)$ instead of $\dheis^m
( L, \nabla; A,A)$ and similarly for the other operator spaces. 
Our goal is to generalize Theorem \ref{theo:resolvent_et_projecteur} for $P
\in \dheis^2 ( L, \nabla; A)$ having a symbol $\si $ of the form 
\begin{gather} \label{eq:symbol} 
 \si (x,\xi ) = \tfrac{1}{2} | \xi |_x ^2 + V(x)
\end{gather}
where $| \cdot |$ is the norm of $T^*M$ for a Riemannian metric of $M$ not
necessarily compatible with $\om$ and $V \in \Ci ( M , \op{End} A)$ is
Hermitian at each point. Example of such operators include Schrödinger
operators with magnetic field and electric potential, holomorphic Laplacians or
{\cc Hodge operators associated to semiclassical Dirac operators}, cf. \cite[Section 3]{oim_copain}. Besides of the numerous examples, the
interest of these operators is that we can compute explicitly the spectrum of
the symbols $\si (x, \cdot)$
$$ \op{sp} ( \si (x, \cdot))  = \Bigl\{ \sum_{i=1}^n B_i (x) ( \al (i) + \tfrac{1}{2}) + V_j (x) / \;
\al \in \N^d, j = 1, \ldots, r \Bigr\} $$
where $0 < B_1 (x) \leqslant \ldots \leqslant B_d (x)$ are the eigenvalues of
$\om(x)$ with respect to $g_x$ and $V_1 (x) \leqslant \ldots \leqslant V_{r} (x)$ are the eigenvalues
of $V(x)$.  Moreover, we have
$$ \tfrac{1}{2} | \xi |^2_x = \sum_{i=1}^d B_i (x) h(s_i, \si_i) , \qquad
h(y,\eta)  = \tfrac{1}{2}  (  y^2 + \eta^2 )$$
where $s_i$ and $\si_i$ are the linear coordinates of $T_x^*M$ associated to a
symplectic basis. So the analysis of $\si ( x, \cdot)$ boils down to the
standard quantum harmonic oscillator $h^w$ or the Landau Hamiltonian $h (
\frac{1}{i} \nabla)$.  

\begin{theo} Let $P \in \dheis^2  (L, \nabla;
  A)$ having a symbol $\si$ of the form \eqref{eq:symbol} and such that for each $k$, $P_k$ is
  formally selfadjoint for a volume element of $M$. Assume $\om$ is nondegenerate and let $\Si = \bigcup_{x \in M} \op{sp} ( \si (x,\cdot)) $.
  Then
  \begin{itemize}
  \item[-]     For any $z \in \C \setminus \Si$, there exists $Q(z) \in \pheis^{-2} ( L,
    \nabla; A)$ such that $Q_k (z) (
    P_k -z ) = \op{id}$ and $ (P_k - z) Q_k ( z) = \op{id}$  when $k$ is large.
  \item[-]  For any $E \in \R \setminus \Si$, $(1_{(-\infty, E]} ( P_k) )$
    belongs to $\pheis^{-\infty} ( L, \nabla ; A)$.  
  \end{itemize}

\end{theo}

The proof is the same as the one of Theorem \ref{theo:resolvent_et_projecteur}. The
symbols $\tau(z)$ and $p_E $ of $Q(z)$ and $1_{(-\infty, E]}$ respectively are
  such that for any $x \in M$, 
$$  \tau(z) (x, \cdot)^w = ( \si (x, \cdot )^w -z )^{-1}, \qquad p_E  (x,\cdot )^w =  1_{(-\infty, E]} ((\si (x, \cdot)^w ) .$$
In the case where $B_i=1$ and $V=0$ , they have been studied for themselves in
\cite{De}, \cite{Un}, and given by the formulas \eqref{eq:Derez} and
\eqref{eq:Un} respectively. {\cc Notice as well that in this case, $\op{sp} (\si
(x) ) = \frac{d}{2} + \N$ and the multiplicity of $\frac{d}{2}
+ m $ is $ { m+d-1 \choose m }$, so that the estimate
\eqref{eq:dim_cluster_asymptotic} follows from Corollary \ref{cor:Weyl_law}.
Moreover, by the first assertion of Theorem \ref{theo:resolvent_et_projecteur}, for any
$\ep>0$ and $M>0$, when $k$ is sufficiently large, 
\begin{gather} \label{eq:clusterweak}
  \op{sp} ( k^{-1} \Delta_k ) \cap (-\infty, M] \subset ( \tfrac{d}{2} + \N
  ) +   (-\ep , \ep ),  
\end{gather}
which is a weak form of \eqref{eq:cluster}.}



\bibliographystyle{plain}
\bibliography{biblio}

\vspace{.5cm}
\noindent
Sorbonne Université, CNRS, IMJ-PRG, F-75006 Paris, France.

\end{document}